\documentclass[11pt, epic, eepic, reqno]{amsart}

 \usepackage{amssymb, amsmath, amsfonts, amsthm, mathrsfs, latexsym, epsfig, psfrag, graphicx}
\usepackage[all]{xy}
\xyoption{v2} \xyoption{knot}\xyoption{arc} \xyoption{poly}
\xyoption{curve}

\newcommand\nc{\newcommand}
\nc\zzzvert {\;|\;}
\nc\zzzcolon {\colon\thinspace}
\nc\zzzsharp {\sharp}
\nc\plim{\varprojlim}
\nc\np{\newpage}
\nc\modone {{\mathbf 1}}
\nc\g{{\mathfrak g}}
\nc\ul{\underline}
\nc\chR{\check{R}}
\nc\rg{\rm range}
\nc\modS {{\mathcal S}}
\nc\modC {{\mathcal C}}
\nc\Ge{\Gamma_k^{\rm even}}
\nc\Gamo{\Gamma_k^{\rm odd}}
\nc\cX{{\mathcal X}}
\nc\ovl{\overline} 
\nc{\dlabel}{\label}
\nc{\dbibitem}{\bibitem}

\newtheorem{theorem}{Theorem}[section]
\newtheorem{lemma}[theorem]{Lemma}
\newtheorem{cor}[theorem]{Corollary}
\newtheorem{proposition}[theorem]{Proposition}
\newtheorem{claim}[theorem]{Claim}

\theoremstyle{definition}
\newtheorem{remark}[theorem]{Remark} 
\theoremstyle{plain}

\newtheorem{definition}[theorem]{Definition}

\newcommand{\Z}{\mathbb Z}

\newcommand{\D}{\Delta}
\nc{\A}{\mathcal A}
\nc{\B}{\mathcal B}

\newcommand{\al}{\alpha}
\newcommand{\be}{\beta}
\newcommand{\gam}{\gamma}
\newcommand{\bbe}{{\ov \beta}}
\newcommand{\bga}{{\ov \gamma}}
\newcommand{\bal}{{\ov \alpha}}

\newcommand{\bxi}{{\ov \xi}}
\newcommand{\R}{\mathbb R}

\newcommand{\N}{\mathbb N}
\nc{\DD}{\mathbb D}
\newcommand{\comment}[1]{}
\newcommand{\spa}{{\rm span}}
\newcommand{\diag}{{\rm diag}}
\newcommand{\ov}{\overline}

\newcommand{\begeq}{\begin{eqnarray*}}
\newcommand{\eneq}{\end{eqnarray*}}
\newcommand{\begal}{\begin{align*}}
\def\ttiny{\@setsize\ttiny{4\p@|\xivpt\@ixvpt}}

 \nc{\scc}{{\text sc}}
 \nc{\oxi}{{\ovl \xi}}
 \nc{\oeta}{{\ovl \eta}}
 \nc{\oal}{{\ovl \al}}
 
\begin{document}

\title[]{Fusion rules on a parametrized series of graphs}
\author{Marta Asaeda}
\thanks{The first named author was sponsored in part by NSF grant  \#DMS-0504199.}

\address{Department of Mathematics, University of California,  Riverside,  900 Big Springs Drive, Riverside, CA, 92521,  USA} 
\email{\tt marta@math.ucr.edu}
\author{Uffe Haagerup}
\address{Department of Mathematical Sciences, University of Copenhagen , 
Universitetspark 5
2100 Copenhagen ¯ Denmark } 
\email{\tt haagerup@math.ku.dk}
\subjclass{46L37}
\keywords{subfactors, fusion algebras}

 \begin{abstract}
 A series of pairs of graphs $(\Gamma_k, \Gamma'_k), k = 0,1,2, \dots $ has been considered as candidates for dual pairs of principal graphs of subfactors of small Jones index above $4$ and it has recently been proved that the pair $(\Gamma_k, \Gamma'_k)$ comes from a subfactor if and only if $k = 0$ or $k =1$. We show that nevertheless there exists a unique fusion system compatible with this pair of graphs for all non-negative integers $k$.
  \end{abstract}
  \maketitle
 \section{Introduction} 
 A  subfactor  $N \subset M$ with finite index and finite depth generates finitely many isomorphism classes of bimodules with four different combinations of left and right coefficients. They form a bi-graded fusion category. Its Grothendieck ring form a {\it fusion ring} or a {\it fusion hypergroup}, namely a bi-graded $\Z$-algebra  ${\mathcal A}$ with following properties: 
 \begin{itemize}
 \item it has a basis given by finitely many irreducible bimodules of four different kinds $\cX = {}_N \cX_N \sqcup {}_N \cX_M \sqcup {}_M \cX_N \sqcup {}_M \cX_M$ (we call the labels $N$, $M$ right or left coefficients, depending on the position), 
 \item an involution $X \in {}_P\cX_Q \to {\ovl X} \in {}_Q \cX_P$ is defined, where $P, Q \in \{N, M\}.$
 \item a product is defined for a pair of bimodules with ``matching" coefficient, namely, for a   pair $(X, Y) \in \cX \times \cX$ such that the right coefficient of $X$ and the left coefficient of $Y$ match, 
 $ XY$ is defined. It decomposes as follows: 
  $$X  Y=\sum N_{X,Y}^Z Z,$$ 
 where the sum is taken over those $Z \in \cX$ that have the same left (resp. right) coefficient as $X$ (resp. $Y$), and  $N_{X,Y}^Z \in \N_0$, moreover Frobenius reciprocity holds: 
 $$N_{X,Y}^Z = N_{Z, {\ovl Y}}^X =N_{{\ovl X}, Z}^Y=N_{{\ovl Y}, {\ovl X}}^{\ovl Z}=N_{{\ovl Z}, X}^{\ovl Y}=N_{Y, {\ovl Z}}^{\ovl X}.$$ 
 \item There are identity objects ${\bf 1}_N \in  {}_N \cX_N$, ${\bf 1}_M \in  {}_M \cX_M$ that act as identity with respect to  the product, whenever it is defined. 
 \end{itemize} 
 The involution extends linearly to define an involution on ${\mathcal A}$. 
 For a fusion ring $\mathcal A$, there is a unique weight function $\mu: {\mathcal A} \to \R_\geq$ satisfying
 \begin{align*}
 \mu({\bf 1}_N)&=\mu({\bf 1}_M)=1, \\ 
 \mu(X  Y)&=\mu(X) \mu (Y), \\
 \end{align*}
 where $X, Y, Z \in \cX$ are with suitable coefficients for each equality, so that $X  Y$ and $X + Z$ are defined. 
  The  {\it  (dual) principal graph} of the subfactor encode partial information of the fusion algebra: namely, the (dual)  principal graph has the vertices corresponding to ${}_N \cX_N \sqcup {}_N \cX_M$ (resp.~ ${}_M \cX_N \sqcup {}_M \cX_M$), with the number of the edges between vertices ${}_N X_N$ and ${}_N Y_M$ (resp.  ${}_M X_M$ and ${}_M Y_N$) given by $N_{X, {}_N M_M}^Y$   (resp. $N_{ X, {}_M M_N}^Y$.) 
 
 On the other hand, one may start with a pair of graphs, and may consider if there is a fusion algebra compatible with the fusion constraints determined by the graphs. Such investigation  may be used to exclude graphs as (dual)  principal graphs of subfactors. For example, type $E_7$ and $D_{2n+1}$ Dynkin diagrams are proved {\it not} to be (dual) principal graphs of subfactors, by showing that the fusion constraints given by the graphs give rise to inconsistency in fusion rules (\cite{I1}, \cite{SV}). Note that the existence of a fusion algebra compatible with a given pair of graphs do not imply the existence of a subfactor with given graphs as (dual) principal graphs. 
 
 In this paper, we deal with the following series of pairs of graphs: 
 
  $$
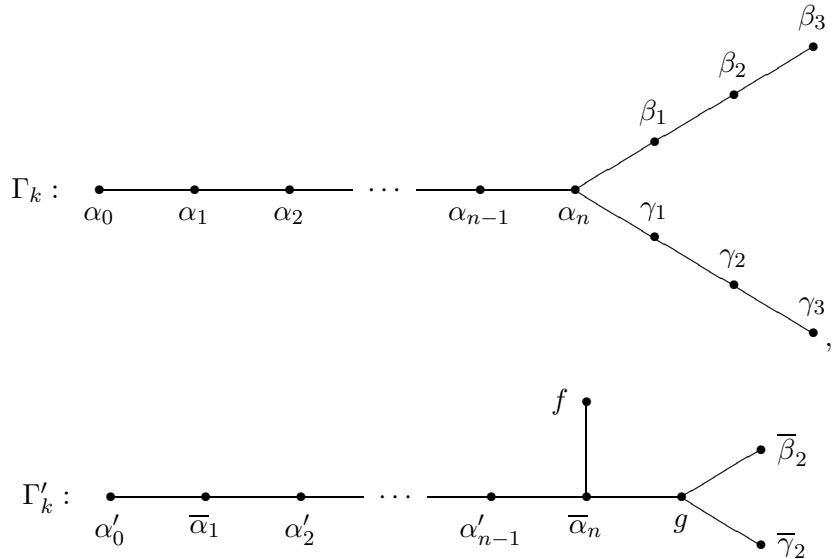
\xy 0;/r1pc/:
 (-1, -5)*{\Gamma_k:};
 (1,-5);(9,-5) **\dir{-};
 (11,-5);(16,-5) **\dir{-};
 (1,-5)*\dir{*}, *++!U{\al_0}; 
 (4,-5)*\dir{*}, *++!U{\al_1}; 
 (7,-5)*\dir{*}, *++!U{\al_2}; 
 (10,-5)*{\cdots};
 (13,-5)*\dir{*}, *++!U{\al_{n-1}}; 
 (16,-5)*\dir{*}, *++!U{\al_n};
  (16,-5);(23.5,-0.5) **\dir{-};
    (16,-5);(23.5,-9.5) **\dir{-};
     (18.5,-3.5)*\dir{*}, *++!D{\beta_1}; (21,-2)*\dir{*}, *++!D{\beta_2}; (23.5,-0.5)*\dir{*}, *++!D{\beta_3}; 
       (18.5,-6.5)*\dir{*}, *++!D{\gam_1}; (21,-8)*\dir{*}, *++!D{\gam_2}; (23.5,-9.5)*\dir{*}, *++!D{\gam_3}; 
\endxy$$
 $$\xy 0;/r1pc/:
 (-1, -5)*{\Gamma_k':};
 (1,-5);(9,-5) **\dir{-};
 (11,-5);(19,-5) **\dir{-};
 (1,-5)*\dir{*}, *++!U{\al'_0}; 
 (4,-5)*\dir{*}, *++!U{{\ovl \al_1}}; 
 (7,-5)*\dir{*}, *++!U{\al'_2}; 
 (10,-5)*{\cdots};
 (13,-5)*\dir{*}, *++!U{\al'_{n-1}}; 
 (16,-5)*\dir{*}, *++!U{{\ovl \al_n}};
  (19,-5)*\dir{*}, *++!U{g};
  (16,-5);(16,-2) **\dir{-};
  (16,-2)*\dir{*}, *++!R{f};
   (19,-5);(21.5,-3.5) **\dir{-};
     (19,-5);(21.5,-6.5) **\dir{-};
     (21.5,-3.5)*\dir{*}, *++!L{{ {\ov \beta_2}}};
         (21.5,-6.5)*\dir{*}, *++!L{{ {\ov \gam_2}}};
\endxy,$$
where $n=4k+3$, $k=0, 1, ...$. Let 
These graphs are a part of the list of the graphs that were candidates for (dual) principal graphs of a subfactor with indices between $4$ and $3+\sqrt{3}$ given by the second author (\cite{H5}). 
Note that the notation used here is somewhat different from the one used in \cite{H5}. It has been already proved that, for $k=0, 1$, $\Gamma_k$ (resp.~ $\Gamma'_k$) are  (dual) principal graphs of a subfactors (\cite{AH}, \cite{BMPS}), and for $k>1$, they are not realized as (dual) principal graphs (\cite{AY}).  
In this paper, we prove that, despite that $\Gamma_k$ (resp.~ $\Gamma'_k$) are not principal graphs for $k>1$, there are still fusion algebras consistent with the graphs, and moreover such fusion algebras are unique for each $k$. Namely we prove the following: 
\begin{theorem}
\dlabel{haupttheorem}
Let 
$V_{11}:=\{\mbox{even vertices of } \Gamma_k\}$, $V_{12}:=\{\mbox{odd vertices of } \Gamma_k\}$, $V_{21}:=\{\mbox{odd vertices of } \Gamma'_k\}$, $V_{22}:=\{\mbox{even vertices of } \Gamma'_k\}$, and $V:=V_{11} \sqcup V_{12} \sqcup V_{21} \sqcup V_{22}$. 
For each $k$, there is a unique fusion algebra ${\mathcal A}=\Z \cX$, where 
$$\cX = {}_N \cX_N \sqcup {}_N \cX_M \sqcup {}_M \cX_N \sqcup {}_M \cX_M$$
compatible with the graphs $\Gamma_k$, $\Gamma'_k$.  Namely
 $$  {}_N \cX_N =V_{11},$$
 $$ {}_N \cX_M  =V_{12},$$
 $${}_M \cX_N=V_{21}, $$
 $$  {}_M \cX_M =V_{22}$$
 as sets, and  
 \begal 
 & N_{X,\al_1}^Y (\mbox{resp. } N_{X, {\bal_1}}^Y) =
 \begin{cases}
 1 \; \mbox{  if $X$ and $Y$ are connected by an edge}  \\
 0 \; \mbox{  else,}
 \end{cases} \\
& N_{X, 1}^Y=\delta_{X,Y}
 \end{align*}
 holds, where $X, Y \in \cX$, and $1$ denotes identity objects $1_N=\al_0 \in {}_N\cX_N$ or $1_M=\al'_0 \in {}_M \cX_M$. 
  \end{theorem}

 The content of this paper is as follows. In Section \ref{unique} we show that if there is a   fusion system compatible with the graphs $\Gamma_k, \Gamma'_k$, it must be unique.   In Section \ref{existense} we show the existence of such a fusion system. 
 
 Part of the work on this paper was conducted while both authors were visiting the Department of Mathematics at Universit\"at M\"unster in April 2009. The authors wish to thank the Department and especially Joachim Cuntz for the hospitality during their stay. 
   
 %
 %
 %
%
%
\section{Uniqueness, positivity, and integrality of the fusion rules}
\dlabel{unique}
In this section we prove that, if there is a fusion algebra compatible with the graphs, it is unique. Positivity and integrality of fusion coefficients is derived: we do not impose them in showing uniqueness of the fusion rules.

 %
\subsection{Fusion rules for the even vertices}
  \dlabel{fusioneven}
 In this subsection we show that there is a unique fusion algebra structure on ${\mathcal A}_1=\Z  {}_N \cX_N$ compatible with the graph $\Gamma_k$. 
 The main issue is to determine the fusion rule among $\be_1, \be_3, \gam_1, \gam_3$. The rest will follow easily from this. 
 
  In the following we assume that there is a fusion algebra compatible with $(\Gamma_k, \Gamma'_k)$.  The involution $\gamma \in V \to {\ovl \gam} \in V$ extends linear to a map on $\R V$. 
 For simplicity, we refer to the objects in $\cX$ by corresponding vertices in $V$.  For $X:= \sum N_X^Z Z \in \R V$ and $Y \in V$, we denote 
 $$<X, Y> = <Y, X> := N_X^Y.$$
 %
 Observe that $< \cdot, \cdot>$ expends linearly to define a bilinear form on $\R V$, and 
 $$<X  Y, Z>=< X , Z  {\ovl Y}>=<Y, {\ovl X}  Z>$$ 
 holds by Frobenius reciprocity. 
The graph $\Gamma_k$ encodes
s
  the decomposition of $X \al_1$ for $X$ in $V_{11}$ into a direct sum of vertices from $V_{12}$ and the decomposition of $Y \bal_1$ into a direct sum of vertices from $V_{11}$. Let $G$ be the adjacency matrix for $(V_{11}, V_{12})$, namely
\begin{align*}
G&=(G_{X,Y})_{X \in V_{11}, Y \in V_{12}},   \\
\mbox {where  } G_{X, Y}&=( \mbox{the number of the edges connecting $X$ and $Y$} ) \\
&=<X  \al_1, Y> \\
&=<Y \bal_1, X>,  
\end{align*}
which is written as the following $(\frac{n+1}{2}+4) \times (\frac{n+1}{2}+2)$-matrix : 
\begin{equation} 
\dlabel{G} 
G=
\bordermatrix{ &  \be_2 & \gam_2  &\al_n    & \al_{n-2}   &\cdots   &  \cdots &  \al_1   \cr
\be_3 &1 & 0 & 0 & 0 & \cdots & \cdots & 0    \cr
\be_1 &1 & 0 & 1 &0& \cdots  & \cdots  & 0&   \cr
\gam_3 & 0 & 1 & 0 & 0& \cdots &   \cdots & 0  \cr
\gam_1 & 0 & 1 & 1 & 0 & \cdots & \cdots & 0    \cr
\al_{n-1}  & \vdots & 0 & 1 & 1 & 0 &  \cdots & 0    \cr
\vdots & \vdots &   \vdots &   \ddots & \ddots & \ddots & \ddots & \vdots    \cr
\al_{2} & 0 & 0& \cdots & 0 & 1 & 1 & 0    \cr
\al_0  & 0 & 0& \cdots & \cdots & 0 & 1 & 1    
  },
  \end{equation}
Let
$${\Delta} :=
\left(\begin{array}{cc}0 & G   \\ G^t & 0 \end{array}\right), 
$$
then 
$$ \Delta^2= 
\left(\begin{array}{cc}G G^t & 0 \\0 & G^t G\end{array}\right).
$$
We put $ {\mathbb D}:=G G^t$, which acts on ${\ovl \A_1}:= \R V_{11}.$ We utilize certain eigen vectors of $\DD$ to determine the fusion structure of $\A_1$. 

Observe from the graph that 
\begin{align*}
& \D \be_1 =\al_n +\be_2, \; \;  \D \gam_1 =\al_n +\gam_2, \\
& \D \be_2 = \be_1+\be_3, \; \; \D \gam_2= \gam_1 +\gam_2, \\
& \D \be_3 =\be_2, \; \; \D \gam_3 =\gam_2. 
\end{align*}
Put 
\begin{align*}
\xi &=(\be_1 -\gam_1) + (\be_3 -\gam_3)\\
\eta & = (\be_1 -\gam_1) - (\be_3 -\gam_3). 
\end{align*}
Then
\begin{align*}
\DD \xi &= \D^2 \xi =\D(2\be_2-2\gam_2)=2\xi, \\
\DD \eta &= \D^2 \eta =0. 
\end{align*}
Let $E(\DD, c),  c\in \R$ be the eigenspace of the eigenvalue $c$ for $\DD$ in $\R(V_{11})$. 
\begin{lemma}
$$\dim E(\DD, 2)=E(\DD, 0)=2
 $$
\end{lemma}
\noindent
{\bf Proof} \\
   
$$
  \DD=
\bordermatrix{ &  \be_3 & \be_1   & \gam_3   & \gam_1      & \al_{n-1} &\cdots&\cdots &  \cdots & \al_{2} &  \al_{0}   \cr
 \be_3    &1 & 1 & 0 & 0 & 0 & 0  & \cdots & \cdots  &\cdots & 0       \cr
 \be_1    &1 & 2 & 0 & 1 & 1  & 0  &   & &  & \vdots       \cr
 \gam_3 & 0 & 0 & 1 & 1 & 0 & 0  &   &  & & \vdots     \cr
 \gam_1  & 0 & 1 & 1 & 2 & 1 & 0 &   &   && \vdots       \cr
\al_{n-1}& 0 & 1 & 0 & 1 & 2 & 1 & 0 &   && \vdots      \cr
\al_{n-3}& 0 & 0 & 0 & 0 & 1 & 2 & 1 & 0  && \vdots      \cr
 \vdots & \vdots &   &   &   & \ddots & \ddots & \ddots & \ddots & \ddots& \vdots      \cr
\vdots & \vdots &   &   &   & & 0 & 1 & 2 & 1 & 0      \cr
 \al_{2} & 0 & \cdots & \cdots & \cdots & \cdots & \cdots  & 0 & 1 & 2 & 1     \cr
 \al_{0} & 0 & \cdots & \cdots & \cdots & \cdots & \cdots & \cdots  & 0 &  1 & 1 } 
$$
Recall that $n=4k+3$. 
Let $\rho_k(x):=\det (tI - \DD)$ be the characteristic polynomial of $\DD=G G^t$. 
It was proved in \cite{A2} that the characteristic polynomial of $G^t G$ is equal to $(t-2)^2 q_k(t)$, where the polynomials $q_k(t), k \geq 0$, can be defined recursively by
 \begeq
q_0(t) &=&t^2-5t+3 \\
q_1(t)&=&   (t-1)(t^3-8t^2+17t-5)  \\
q_k(t)&=&(t^2-4t+2)q_{k-1}(t)-q_{k-2}(t), \; k \geq 2 \\
\eneq
Since the matrix $G$ has $2k+6$ rows and $2k+4$ columns, $G G^t$ is a unitary conjugate of $G^t G \oplus O_2$, where $O_2$ is the zero $2 \times 2$ matrix. Hence 
\begin{align*}
\dlabel{formularho} 
\rho_k(t)&=t^2 \det(tI-G^tG) \\
&=t^2 (t-2)^2 q_k(t). 
\end{align*}
Using the recursion formula for $q_k(t)$, one gets
$Êq_k(0) = 2k+3$ Êand $q_k(2) = (-1)^(k+1) (2k+3)$
 In particular neither $0$ nor $2$ is a root of $q_k$. Hence $0$ and $2$ are roots of multiplity $2$ in $\rho_k$. Since $\DD = G G^t$
Êis a symmetric matrix, it follows that the dimensions of the eigenspaces for $\DD$ for the eigenvalues $0$ and $2$ are both equal to two.
%


Bases of $E(\DD, 2)$, $E(\DD, 0)$ may be taken as follows: 
\begin{align*}
E(\DD, 2):=\spa\{x_1, x_2\} \\
E(\DD, 0):=\spa\{y_1, y_2\}, 
\end{align*}
where
\begin{align*}
x_1& :=2(\al_0 +\al_2)-2 (\al_4 +\al_6)+ \cdots + (-1)^k 2(\al_{4k} +\al_{4k+2})  \\
& + (-1)^{k+1} (\be_1 +\gam_1 +\be_3 +\gam_3) \\
x_2& := \xi =(\be_1-\gam_1) +(\be_3-\gam_3) \\
y_1&:= 2\al_0 -2 \al_2 + \cdots + 2\al_{4k}-2\al_{4k+2} + (\be_1+\gam_1)  - (\be_3 +\gam_3)\\
y_2 &:= \eta =(\be_1 -\gam_1) -(\be_3 -\gam_3)
\end{align*}
Assume that we have a fusion algebra compatible with the pair of the graphs $(\Gamma_k, \Gamma'_k)$, 
and let $\pi$ and $\pi'$ be the conjugate maps $\gamma \mapsto {\ovl \gamma}$  on $V_{11}$ and $V_{22}$. By the argument used in \cite[pp28-31]{H5}, $\pi'$ fixes every element of $V_{22}$. For $\pi$, there are only two possibilities: 
\begin{itemize}
\item[{\bf Case 1}] (=case (b) in \cite[p31]{H5}) 
$$ {\ovl \be_1}=\be_1, {\ovl \gam_1}=\gam_1, {\ovl \be_3}=\gam_3 ( \Leftrightarrow {\ovl \gam_3}=\be_3), $$
\item[{\bf Case 2}] (=case (a) in \cite[p31]{H5}. To be eliminated.) 
$$ {\ovl \be_1}=\gam_1  ( \Leftrightarrow {\ovl \gam_1}=\be_1), {\ovl \be_3}=\be_3, {\ovl \gam_3}=\gam_3.$$
\end{itemize}
In both cases, ${\ovl \al_{2j}}=\al_{2j}$ for $ j=0, 1, \ldots, 2k+1$. 
Note that $\pi$ extends linearly to ${\mathcal A}_1$ and ${\ovl {\mathcal A}_1}=\R V_{11}$. 
Let $E(\DD, c)_{{\text sc}}:=   E(\DD, c)^\pi$. Observe that 
$$ c_1 {\ovl x_1} +c_2 {\ovl x_2} =c_1 x_1+ c_2 x_2, \; c_1, c_2 \in \R$$
holds if and only if $c_2=0$ in both cases 1 and 2, and similarly
$$c_1 c_1 {\ovl y_1} +c_2 {\ovl y_2} =c_1 y_1+ c_2 y_2, \; c_1, c_2 \in \R$$
if and only if $c_2=0$ in both cases. Therefore
\begin{align*}
E(\DD, 2)_{{\text sc}}&=\R x_1 \\
E(\DD, 0)_{{\text sc}}&=\R y_1
\end{align*}
By the definition of principal graphs,  the matrix  $\DD: \R V_{11} \to \R V_{11}$ corresponds to the fusion rule of  the right tensor product by $\al {\ovl \al}$, where $\al=\al_1$. Therefore
\begin{align*}
\DD({\ovl \xi} \xi)&= {\ovl \xi} \DD(\xi) =2 {\ovl \xi} \xi \\
 \DD({\ovl \eta} \eta)&= {\ovl \eta} \DD(\eta) =0.
\end{align*}
  Hence 
\begin{align*}
{\ovl \xi} \xi \in E(\DD, 2)_{{\text sc}} &=\R x_1, \\
{\ovl \eta} \eta \in E(\DD, 0)_\scc &= \R y_1.
\end{align*}
 Thus 
 $$<{\ovl \xi}\xi, \al_0>=<\xi, \xi \al_0> =<\xi, \xi> =4.$$
 Hense the coefficient of ${\ovl \xi} \xi$ at $\al_0$ is 4. Since $\oxi \xi \in \R x_1$, we have $\oxi \xi =2 x_1$. Likewise we obtain $\oeta \eta =2y_1.$ Noting that 
 $$\oxi = 
 \begin{cases} 
 \eta \; \mbox{ for Case 1} \\
- \eta \; \mbox{ for Case 2} , 
 \end{cases}
 $$
 we have 
 \begin{itemize}
 \item \mbox{In Case 1:} $\xi \eta =2y_1, \eta \xi=2 x_1, $\\
 \item \mbox{In Case 2:} $\xi \eta =-2y_1, \eta \xi=-2 x_1$.
 \end{itemize}
 \begin{lemma}
 $$\xi^2=0,  \eta^2=0.$$
 \end{lemma}
 {\bf Proof}. 
Since $\DD(\xi^2)=\xi \DD(\xi)=2\xi^2$, $\xi^2=c_1 x_1 +c_2 x_2$ for some $c_1, c_2 \in \R.$ Moreover, since $<\xi, \eta>=0$, we have
 \begin{align*}
 <\xi^2, \al_0> = < \xi, \oxi \al_0> &=\pm <\xi, \eta> \\
 &=0
 \end{align*}
 Together with $<c_1 x_1 + c_2 x_2, \al_0>=2c_1$, $c_1, c_2 \in \R$, we obtain
 $$\xi^2 =c_2 x_2=c_2 \xi.$$
 We show that $c_2=0$:
 \begin{align*}
4 c_2&=<c_2 \xi, c_2 \xi>= <\xi^2, \xi^2> =<\oxi \xi, \xi \oxi>=4<x_1, y_1> \\
&=(2-2)-(2-2)+ \cdots (-1)^k(2-2)+(1+1-1-1)=0.
\end{align*}
We used that $\oxi \xi=2x_1$, $\xi \oxi=2 y_1$ for both cases. Thus $\xi^2=0$. Then $\oxi^2=\eta^2=0$ for both cases. \qed

Since $\be_3-\gam_3 =\frac{1}{2}(\xi-\eta)$, we get 
\begin{align*}
(\be_3-\gam_3)^2 &= \frac{1}{4}(\xi-\eta)^2  \\
&= \frac{1}{4} (\xi^2 + \eta^2 -\xi \eta -\eta \xi) \\
&= -\frac{1}{4} (\xi \eta +\eta \xi)  \\
&= 
\begin{cases}
-\frac{1}{2} (x_1 +y_1) \mbox{  in Case 1} \\
\frac{1}{2} (x_1 +y_1)  \mbox{  in Case 2} \\
\end{cases}
\end{align*}
\begin{remark}
\dlabel{x1+y1}
For $k$ even (i.e. $n=3$ (mod $8$)) and $k=2l$, 
 $$
\frac{1}{2} (x_1 +y_1)  
 = 
 2( \al_0 -\al_6+\al_8 -\al_{14} +\al_{16} - \cdots + \al_{8l}) - (\be_3 +\gam_3) $$
 and  for $k$ odd (i.e. $n=7$ (mod $8$)) and $k=2l+1$, 
\begal
&\frac{1}{2} (x_1 +y_1)  \\
&= 2( \al_0 -\al_6+\al_8 -\al_{14} +\al_{16} - \cdots + \al_{8l} -\al_{8l+6}) + (\be_1 +\gam_1)
\end{align*}
 \end{remark}
Consider next the sequence of polynomials $R_n$ given recursively by 
\begeq
\dlabel{r}
\notag R_0(t)& =1, \\
\notag R_1(t)& =t, \\
R_m(t)&=tR_{m-1}(t)-R_{m-2}(t), n \geq 2. 
\eneq
as in \cite[p33--34]{H5}. 
Note that $R_m(t)=U_m(\frac{t}{2})$, where $U_m$ is the $m$-th Chebyshev polynomial of second kind \cite{ch}. Moreover, 
$$R_m(2 \cos \theta)=\frac{\sin (m+1)\theta}{\sin \theta}, \; 0 < \theta < \pi. $$
 By the recursion formula for $R_n$, it follows that 
 \begin{align*}
 R_j(\D) \al_0&=\al_j, \; 0 \leq j \leq n, \\
 R_{n+1}(\D) \al_0&= \be_1 +\gam_1, \\
 R_{n+2}(\D) \al_0&= \al_n +\be_2 +\gam_2, \\
 R_{n+3}(\D) \al_0 &= \al_{n-1} +\be_1 +\gam_1 +\be_3 +\gam_3. 
 \end{align*}
 Hence
 \begal 
 \be_3 +\gam_3 &= (R_{n+3}(\D) -R_{n+1}(\D)-R_{n-1} (\D) ) \al_0 \\
 &= (R_{4k+6}(\D) -R_{4k+4}(\D)-R_{4k+2} (\D) ) \al_0
\end{align*}
For $m$ even, $R_m(t)$ is an even polynomial in $t$, thus there is are unique polynomials
$(Q_j)_{j=0,1,2, \ldots} $ with $\deg (Q_l)=l$, such that 
$$Q_j(t^2)=R_{2j}(t), \; t \in \R,  \; j=0,1,2, \ldots. $$
With this notation, we have 
\begeq
\dlabel{b3+g3} \\
 \notag \be_3+\gam_3 &=  (Q_{2k+3}(\DD) -Q_{2k+2}(\DD)-Q_{2k+1} (\DD) ) \al_0 \\
&=(Q_{2k+3}-Q_{2k+2}-Q_{2k+1})(\al \oal).   
\eneq
Therefore 
\begal
(\be_3-\gam_3)(\be_3+\gam_3) &  = (Q_{2k+3}-Q_{2k+2}-Q_{2k+1})(\DD) (\be_3-\gam_3)    \\
&= \frac{1}{2} (Q_{2k+3}-Q_{2k+2}-Q_{2k+1})(\DD)(\xi-\eta)
\end{align*}
Since $\DD \xi =2\xi$ and 
\begal
Q_m(2)=R_{2j}(\sqrt{2})&=\frac{\sin (2j+1)\pi/4 }{\sin \pi /4} \\
&= \begin{cases}
1, \; \;  j=0, 1 \; (\mbox{mod } 4) \\
-1,  \; j=2, 3  \; (\mbox{mod } 4), 
\end{cases}
\end{align*}
we have 
\begal
Q_j(\DD) \xi = 
\begin{cases}
\xi,  \; \;  j=0, 1 \; (\mbox{mod } 4) \\
-\xi,   \; j=2, 3  \; (\mbox{mod } 4), 
\end{cases}
\end{align*}
Similarly, since $\DD \eta=0$ and 
$$Q_j(0)=R_{2j}(0)=\frac{\sin (2j+1) \pi /2}{\sin \pi /2} =(-1)^j,$$
we have 
$$Q_j(\DD)\eta =(-1)^j \eta, \; j=0, 1, 2 \ldots.$$
Therefore we have
\begal
& (Q_{2k+3}(\DD) -Q_{2k+2}(\DD)-Q_{2k+1} (\DD) ) \xi  \\
&=\begin{cases} 
(Q_{4l+3}(\DD) -Q_{4l+2}(\DD)-Q_{4l+1} (\DD) )  \xi = -\xi 
\mbox{                  for $k=2l$, $l \in \N_0$}  \\
(Q_{4l+5}(\DD) -Q_{4l+4}(\DD)-Q_{4l+3} (\DD) )  \xi = \xi  
\mbox{                  for   $k=2l+1$, $l \in \N_0$,}  \\
\end{cases}
\end{align*}
and in both cases
$$ (Q_{2k+3}(\DD) -Q_{2k+2}(\DD)-Q_{2k+1} (\DD) ) \eta=-\eta.$$ 
Hense 
\begal
&(\be_3-\gam_3)(\be_3+\gam_3) = \frac{1}{2} (Q_{2k+3}-Q_{2k+2}-Q_{2k+1})(\DD)(\xi-\eta) \\
&= 
\begin{cases}
\frac{1}{2}(-\xi+\eta)=\gam_3 -\be_3, \; k \mbox{ even}, \\
\frac{1}{2}(\xi+\eta)= \be_1 -\gam_1,  \; k \mbox{ odd}. 
\end{cases}
\end{align*}
Using the contragradient map we get \\
\noindent
{\bf For Case 1}:\\
\begal
(\be_3+\gam_3) (\be_3-\gam_3)&= {\ovl {(\bbe_3-\bga_3)(\bbe_3+\bga_3)}} \\
&={\ovl {(\gam_3 -\be_3)(\gam_3 + \be_3)}} \\
&= -{\ovl {(\be_3 -\gam_3)(\be_3+\gam_3)}} \\
&= \begin{cases}
-(\bga_3-\bbe_3)=-(\be_3-\gam_3), \; k \mbox{ even}, \\
-(\bbe_1 -\bga_1)=-(\be_1 -\gam_1), \; k \mbox{ odd},   
\end{cases}
\end{align*}
\noindent
{\bf For Case 2} (to be eliminated): \\
\begal
(\be_3+\gam_3) (\be_3-\gam_3)&= {\ovl {(\bbe_3-\bga_3)(\bbe_3+\bga_3)}} \\
&={\ovl {(\be_3 -\gam_3)(\be_3 + \gam_3)}} \\
\begin{cases}
\bga_3-\bbe_3=\gam_3-\be_3, \; k \mbox{ even}, \\
\bbe_1 -\bga_1=\gam_1 -\be_1, \; k \mbox{ odd}. 
\end{cases}
\end{align*}
Thus in both cases 
$$(\be_3+\gam_3) (\be_3-\gam_3) =
\begin{cases}
\gam_3-\be_3, \; k \mbox{ even}, \\
\gam_1 -\be_1, \; k \mbox{ odd}. 
\end{cases}
$$
So far, we have obtained the following three formulae: 
\begin{itemize}
\item[{\bf [A]}] 
$$(\be_3-\gam_3)^2=
\begin{cases}
-\frac{1}{2} (x_1-y_1) \; \mbox{in Case 1}\\
\frac{1}{2} (x_1-y_1) \; \mbox{in Case 2}
\end{cases}
$$
\item[{\bf [B]}]
$$ (\be_3-\gam_3)(\be_3+\gam_3)  
= 
\begin{cases}
\frac{1}{2}(-\xi+\eta)=\gam_3 -\be_3, \; k \mbox{ even}, \\
\frac{1}{2}(\xi+\eta)= \be_1 -\gam_1,  \; k \mbox{ odd}. 
\end{cases}
 $$
 \item[{\bf [C]}]
 $$(\be_3+\gam_3) (\be_3-\gam_3) =
\begin{cases}
\gam_3-\be_3, \; k \mbox{ even}, \\
\gam_1 -\be_1, \; k \mbox{ odd}. 
\end{cases}
$$
\end{itemize}
Next we compute $(\be_3+\gam_3)^2$, in order to find $\be_3^2$, $\gam_3^2$, $\be_3 \gam_3$ and $\gam_3 \be_3$. 
\begin{claim}  \dlabel{D} \ \\
\begin{itemize}
\item[{\bf [D]}]
$(\be_3 +\gam_3)^2=2(c_0 \al_0 +c_1 \al_2 + \cdots + c_{2k+1} \al_{4k+2}) +c_{2k+2}(\be_1 + \gam_1) 
+c_{2k}(\be_3 +\gam_3)$, 
\end{itemize}
where $c_j$'s are defined by $c_0=1$, $c_1 =c_2=0$ and $c_j=c_{j-1}+c_{j-2}+c_{j-3}$ for $j \geq 3$. 
\end{claim}
{\bf Proof}\\
Recall that 
\begal
(\be_3+\gam_3)&=(Q_{2k+3}-Q_{2k+2}-Q_{2k+1})(\DD) \al_0 \\
&= (R_{4k+6}(\D) -R_{4k+4}(\D)-R_{4k+2} (\D) ) \al_0, 
\end{align*}
thus 
$$(\be_3+\gam_3)^2=(R_{4k+6}(\D) -R_{4k+4}(\D)-R_{4k+2} (\D) ) (\be_3+\gam_3).  \; \; \cdots (\sharp)$$
Our strategy of the proof is as follows: first we find a sequence of polynomials $(S_{j})$ such that $S_j(\D)(\be_3+\gam_3)$ is given by a simple formula. Next we rewrite the right hand side of $(\sharp)$ using $(S_j)$'s. 

Observe that  we obtain from the graph the following: 
\begal
R_0(\D)(\be_3+\gam_3)&=(\be_3+\gam_3), \\
R_1(\D)(\be_3+\gam_3)&=(\be_2+\gam_2), \\
R_2(\D)(\be_3+\gam_3)&=\D (\be_2+\gam_2) -(\be_3+\gam_3) =\be_1+\gam_1, \\
R_3(\D)(\be_3+\gam_3)&=\D (\be_1+\gam_1) -(\be_2+\gam_2) = 2\al_n, \\
R_4(\D) (\be_3+\gam_3)&=2 \D \al_n -(\be_1+\gam_1)=2 \al_{n-1} +\be_1+\gam_1,  
\end{align*}
We define the polynomials $(S_{j}(t) )_{j\geq 3}$ by the following recursive formula: 
\begin{align*}
S_3(t)&=R_3(t),  \\
S_4(t)&=R_4(t)-R_2(t), \\
S_j(t)&=t S_{j-1}(t) - S_{j-2}(t), \; j \geq 5. 
\end{align*}
By definition $S_3(\D) (\be_3+\gam_3)=2\al_n$, $S_4(\D)(\be_3+\gam_3)=2 \al_{n-1}$. 
Since $\al_{l-1} = \D\al_l -\al_{l+1}$ for $l=1, 2, \ldots, n-1$, we easily obtain 
$$S_j(\D)(\be_3+\gam_3)=2 \al_{n-j+3}$$ 
for $j=3, 4, \ldots, n+3.$
Next we express $R_j$'s in terms of $S_j$'s. 
\begin{lemma}
\dlabel{SRrelations}
For $j \geq 2$, 
\begal
R_{2j-1}&= d_0 S_{2j-1} +d_1 S_{2j-3} + \cdots + d_{j-2} S_3 + (d_{j-1} -d_{j-2})R_1 \\
R_{2j}&=d_0 S_{2j} +d_1 S_{2j-2} + \cdots + d_{j-2} S_4 + d_{j-1} R_2 +d_{j-3}R_0, 
\end{align*}
where $d_j$'s satisfy 
\begal
d_j&=d_{j-1}+d_{j-2} + d_{j-3} \\
d_{-1}&=0, d_0=d_1=1,  
\end{align*}
\end{lemma}
{\bf Proof of Lemma}:\\
For $j=2$ it is obvious by the definition of $S_j$'s. We proceed with induction. Assume that it is true for $j$ ($j \geq 2$). Using the recursion formulae for $R_j$'s and $S_j$'s, we have
\begal
&R_{2j+1}(t)= t R_{2j}(t)-R_{2j-1}(t) \\
&= t (d_0 S_{2j} +d_1 S_{2j-2} + \cdots + d_{j-2} S_4 + d_{j-1} R_2 +d_{j-3}) \\
& \hspace{0.5cm}- ( d_0 S_{2j-1} +d_1 S_{2j-3} + \cdots + d_{j-2} S_3 + (d_{j-1} -d_{j-2})R_1)  \\
&= d_0 S_{2j+1} + d_1 S_{2j-1} + \cdots d_{j-2} S_5 + t (d_{j-1}R_2 +d_{j-3}) - (d_{j-1} -d_{j-2})R_1 \\
&=d_0 S_{2j+1} + d_1 S_{2j-1} + \cdots d_{j-2} S_5 + d_{j-1} (t R_2 - R_1) + t d_{j-3} - d_{j-2}R_1 \\
&= d_0 S_{2j+1} + d_1 S_{2j-1} + \cdots d_{j-2} S_5 + d_{j-1} S_3 +(d_{j-3}-d_{j-2})R_1. 
\end{align*}
The last equality was obtained using $S_3=R_3$, $R_1=t$,  and  $d_{j-2} +d_{j-3}=d_j -d_{j-1}$. 
Likewise we have
\begal
R_{2j+2}(t)& =tR_{2j+1}(t)-R_{2j}(t) \\
&= d_0 S_{2j+2} + d_1 S_{2j} + \cdots  + d_{j-2} S_6 \\
& \hspace{1cm} + t( d_{j-1} S_3 +(d_{j}-d_{j-1})R_1)  - (d_{j-1} R_2 +d_{j-3}R_0)  \\
&=d_0 S_{2j+2} + d_1 S_{2j} + \cdots  + d_{j-2} S_6  + d_{j-1} R_4   \\
& \hspace{1cm}  + (d_{j}-d_{j-1}) (R_2+R_0)-d_{j-3}R_0 \\
&=d_0 S_{2j+2} + d_1 S_{2j} + \cdots  + d_{j-2} S_6  + d_{j-1} S_4  +  \\
&  \hspace{1cm} +  d_j R_2 +   (d_{j}-d_{j-1}- d_{j-3})  R_0 \\
 &=d_0 S_{2j+2} + d_1 S_{2j} + \cdots  + d_{j-2} S_6  + d_{j-1} S_4    +  d_j R_2 + d_{j-2} R_0.
\end{align*}
\qed

Let us go back to ($\sharp$). Using Lemma \ref{SRrelations}, 
\begal 
& R_{4k+6}  -R_{4k+4} -R_{4k+2}   \\
 & = d_0 S_{4k+6} + (d_1-d_0) S_{4k+4} + d_{-1} S_{4k+2} + d_0 S_{4k} + d_1 S_{4k-2} + \cdots  \\
 &\hspace{1cm}  + d_{2k-2} S_4 + d_{2k-1} R_2 + d_{2k-3} R_0\\
&=  S_{4k+6}  +    d_0 S_{4k} + d_1 S_{4k-2} + \cdots  \\
 &  \hspace{1cm} + d_{2k-2} S_4 + d_{2k-1} R_2 + d_{2k-3} R_0
 \end{align*} 
 Recall 
 \begal
 S_j(\D) (\be_3 +\gam_3) &=2 \al_{n-j+3} \\
 R_2 (\be_3+\gam_3)  & =\be_1 +\gam_1. 
 \end{align*}
 Letting  $c_0:=1$, $c_1=c_2=0$, $c_j:=d_{j-3}$ for $j \geq 3$, we obtain the formula {\bf [D]}. This concludes the proof for Claim \ref{D}.  \qed
 
Thus far we obtained the formulae for $(\be_3 - \gam_3)^2$, $(\be_3-\gam_3)(\be_3+\gam_3)$, $(\be_3+\gam_3)(\be_3-\gam_3)$ and $(\be_3+\gam_3)^2$ as in {\bf [A]},  {\bf [B]},  {\bf [C]},  {\bf [D]}. This enable us to understand the fusion rules among $\be_3, \gam_3$ and their conjugates. We obtain the following: 
\begin{proposition}
The Case 2 does not occur. Namely $\be_1$, $\gam_1$ are self conjugate and $\bbe_3=\gam_3$ if there is a fusion algebra compatible with the graphs $\Gamma_k$,  $\Gamma'_k$. 
\end{proposition}
\noindent
{\bf Proof.}\\
First observe that, by the definition $(c_j)_{j \geq 0}$ used in Claim \ref{D}, it follows that $c_j$ (mod $4$) is periodic in $j$ with period $8$. The values are given in the following Table \ref{cjmod4}:  \\
\begin{table}[h]
\caption{}
\begin{tabular}{|c|c|c|c|c|c|c|c|c|}
\hline $j$ (mod $8$)  & 0 & 1 & 2 & 3 & 4 & 5 & 6 & 7 \\
\hline $c_j$ (mod $4$) & 1 & 0 & 0 & 1 & 1 & 2 & 0 & 0 \\
\hline \end{tabular}
\dlabel{cjmod4}
\end{table}

 In particular, 
\begin{equation*}
\dlabel{stars}
(\star)
\begin{cases}
c_{2j}=1 \; (\mbox{mod } 4) \mbox{ for } $j$ \mbox{ even}, \\
c_{2j}=0 \; (\mbox{mod } 4) \mbox{ for } $j$ \mbox{ odd}. \\
\end{cases}
\end{equation*}
In the following we assume Case 2 and derive contradiction.  \\
\noindent
$\bullet$ {\bf for $k$ even:} 

By {\bf [B]} and {\bf [C]}, we have 
$$  (\be_3-\gam_3)(\be_3+\gam_3)=(\be_3+\gam_3)(\be_3-\gam_3),$$
hence 
\begal
\be_3 \gam_3 =\gam_3 \be_3 & =\frac{1}{2} (\be_3\gam_3 +\gam_3 \be_3) \\
&=\frac{1}{4}((\be_3 + \gam_3)^2 -(\be_3 - \gam_3)^2). 
\end{align*}
From {\bf [A]} (Case 2),  {\bf [D]} and Remark \ref{x1+y1}, it follows that the coefficient of $\be_3$ in the expansion of $\be_3 \gam_3$ in irreducible objects is equal to 
$$\frac{c_{2k} +1}{4}.$$
 Since $k$ is even, $c_{2k}=1$ mod $4$ by $(\star)$, $(c_{2k} +1)/4$ is not an integer. This implies that Case 2 does not occur if $k$ is even. 
 
 \noindent
 $\bullet$ {\bf for $k$ odd}: \\
 From {\bf [B]}, {\bf [C]}, we get 
 $$  (\be_3-\gam_3)(\be_3+\gam_3)=- (\be_3+\gam_3)(\be_3-\gam_3).$$
 Hence 
 \begal
 \be_3^2 =\gam_3^2 &= \frac{1}{2} (\be_3^2 +\gam_3^2)\\
 &=\frac{1}{4} ((\be_3 + \gam_3)^2 +(\be_3 - \gam_3)^2).
 \end{align*}
 From {\bf [A]} (Case 2),  {\bf [D]} and Remark \ref{x1+y1}, it follows that the coefficient of $\be_1$ in the expansion of $\be_3^2 $ in irreducible objects is equal to 
$$\frac{c_{2k+2} +1}{4}.$$
 Since $k$ is odd, $c_{2k+2}=1$ mod $4$ by $(\star)$, $(c_{2k} +1)/4$ is not an integer.  This excludes Case 2 for $k$ odd as well. \qed
 
 In the following we determine all the irreducible decompositions for the products of any two objects in $V$, and show that the coefficients are non-negative integers. Since we excluded Case 2, we rewrite the formula {\bf{[A]}}:
 \begin{itemize}
 \item[{\bf{[A']}}] 
 For $k=2l$, $l=0, 1, 2, \ldots,$
 $$(\be_3-\gam_3)^2= -2( \al_0 -\al_6+\al_8 -\al_{14} +\al_{16} - \cdots + \al_{8l}) - (\be_3 +\gam_3), $$
 and for $k=2l+1$, $l=0, 1, 2, \ldots,$
 $$(\be_3-\gam_3)^2= -2 ( \al_0 -\al_6+\al_8 -\al_{14} +\al_{16} - \cdots + \al_{8l} -\al_{8l+6}) + (\be_1 +\gam_1).$$
 \end{itemize}
Put 
\begal
A&: = (\be_3-\gam_3)^2 \\
B&: = (\be_3-\gam_3)(\be_3+\gam_3) \\
C&:= (\be_3+\gam_3)(\be_3-\gam_3) \\
D&:= (\be_3+\gam_3)^2.
\end{align*}
Then 
\begal
\be_3 \gam_3 &= \frac{(D-A) + (B-C)}{4} \\
\gam_3 \be_3 &=  \frac{(D-A) - (B-C)}{4}  \\
\be_3^2&=  \frac{(D+A) + (B+C)}{4} \\
\gam_3^2&= \frac{(D+A) - (B+C)}{4} 
\end{align*}
 We introduce new constants $(f_j)_{j \geq 0}$, $(g_j)_{j \geq 0}$ by 
 \begal
 \begin{cases}
 f_j=\frac{1}{2} (c_j+1), g_j=\frac{1}{2} (c_j-1) \mbox{ when } j=0 (\mbox{mod } 4),\\
  f_j=\frac{1}{2} (c_j-1), g_j=\frac{1}{2} (c_j+1) \mbox{ when } j=3 (\mbox{mod } 4),\\
   f_j= g_j=\frac{1}{2} c_j  \mbox{ when } j=1, 2 (\mbox{mod } 4).\\
   \end{cases}
   \end{align*}
   Note that  $f_j + g_j=c_j$ for all $j$. Furthermore, from Table \ref{cjmod4}, observe that $f_j$, $g_j$'s are non-negative integers for all $j \geq 0$. The list of some values for $f_j$'s and $g_j$'s are given in Table \ref{fjgj}:
\begin{table}[h]
\caption{}
\begin{tabular}{|c|c|c|c|c|c|c|c|c|c|c|c|c|c|}
\hline $j$   & 0 & 1 & 2 & 3 & 4 & 5 & 6 & 7   & 8 & 9 & 10 & 11 &12  \\
\hline $f_j$ & 1 & 0 & 0 & 0 & 1 & 1 & 2 & 3 &7 & 12 & 22 & 40 & 75  \\
\hline $g_j$ & 0 &0 &0 & 1 & 0 & 1 & 2 & 4 &6 & 12 & 22 &41 & 74  \\
\hline \end{tabular} 
\dlabel{fjgj}
\end{table}

For $k$ even, using the formulae {\bf{[A']}},  {\bf{[B]}},  {\bf{[C]}},  {\bf{[D]}}, we have 
\begal
\frac{D-A}{4}&= f_0 \al_0 + f_1 \al_2 + \cdots f_{2k+1}\al_{4k+2}  \\
& \hspace{1cm} +\frac{1}{4} c_{2k+2} (\be_1 +\gam_1) +\frac{1}{4} (c_{2k}-1)(\be_3 +\gam_3), \\
\frac{D+A}{4} &= g_0 \al_0 + g_1 \al_2 + \cdots + g_{2k+1} \al_{4k+2}  \\
& \hspace{1cm} +  \frac{1}{4} c_{2k+2} (\be_1 +\gam_1) + \frac{1}{4} (c_{2k}+1) (\be_3 +\gam_3), \\
\frac{B-C}{4}& =0,  \\
\frac{B+C}{4}&= \frac{1}{2}(\gam_3 -\be_3). 
\end{align*}
Since $k$ is even, 
$c_{2k+2} = 2 f_{2k+2} =2 g_{2k+2}$, $c_{2k}+1 =2f_{2k}$, $c_{2k} -1 =2g_{2k}$. Hence we obtain the following theorem:
\begin{theorem}
\dlabel{b3g3fusionEVEN} 
For $k$ even, 
\begal
\be_3 \gam_3 & =\gam_3 \be_3=  f_0 \al_0 + f_1 \al_2 + \cdots f_{2k+1}\al_{4k+2}  \\
& \hspace{2cm} +\frac{1}{2} f_{2k+2} (\be_1 +\gam_1) +\frac{1}{2} (f_{2k}-1)(\be_3 +\gam_3), \\
\be_3^2 &= g_0 \al_0 + g_1 \al_2 + \cdots + g_{2k+1} \al_{4k+2}  \\
& \hspace{2cm} + \frac{1}{2} g_{2k+2} (\be_1 +\gam_1) + \frac{1}{2} g_{2k} \be_3 + \frac{1}{2} (g_{2k}+2) \gam_3, \\
\gam_3^2 &= g_0 \al_0 + g_1 \al_2 + \cdots + g_{2k+1} \al_{2k+2}  \\
& \hspace{2cm} +  \frac{1}{2} g_{2k+2} (\be_1 +\gam_1) + \frac{1}{2} (g_{2k}+2) \be_3 + \frac{1}{2} g_{2k} \gam_3.
\end{align*}
All the coefficients of irreducible elements are non-negative integers. 
\end{theorem}
{\bf Proof}. The only remaining thing to prove is that $f_{2k+2}$ is even, $f_{2k}$ is odd, $g_{2j}$ is even for any $j$. 
Since $k$ is even, $c_{2k+2} =0$ (mod $4$). Thus
$f_{2k+2}= \frac{1}{2} c_{2k+2}$ is even.  Likewise $c_{2k}=1$ (mod $4$), thus $f_{2k} =\frac{1}{2} (c_{2k} +1)$ is odd. 
$$g_{2j}=
\begin{cases}
\frac{1}{2} (c_{2j} -1)  \mbox{ for } j \mbox{ even}, \\
\frac{1}{2} c_{2j}   \mbox{ for } j \mbox{ odd}
\end{cases}
$$
Since $c_{2j}-1=0$ (mod $4$) for $j$ even, $c_{2j}=0$ (mod $4$) for $j$ odd, $g_{2j}$ is even for any $j$. \qed
 
In the same way, we get  for $k$ odd: 
\begal
\frac{D-A}{4}&= f_0 \al_0 + f_1 \al_2 + \cdots f_{2k+1}\al_{4k+2}  \\
& \hspace{1cm} +\frac{1}{4} (c_{2k+2}+1) (\be_1 +\gam_1) +\frac{1}{4} c_{2k} (\be_3 +\gam_3), \\
\frac{D+A}{4} &= g_0 \al_0 + g_1 \al_2 + \cdots + g_{2k+1} \al_{2k+2}  \\
& \hspace{1cm} +  \frac{1}{4} (c_{2k+2}-1) (\be_1 +\gam_1) + \frac{1}{4} c_{2k} (\be_3 +\gam_3), \\
\frac{B-C}{4}& =\frac{1}{2}(\be_1 -\gam_1),  \\
\frac{B+C}{4}&= 0. 
\end{align*}
Since $k$ is odd, $c_{2k+2}+1 = 2 f_{2k+2}$, $c_{2k+2}-1 =2 g_{2k+2}$, $c_{2k}  =2f_{2k} =2g_{2k}$. Hence we get: 
\begin{theorem}
\dlabel{b3g3fusionODD}
For $k$ odd, 
\begal
\be_3 \gam_3 & =  f_0 \al_0 + f_1 \al_2 + \cdots f_{2k+1}\al_{4k+2}  \\
& \hspace{2cm} +\frac{1}{2}( f_{2k+2}+1) \be_1 +\frac{1}{2}( f_{2k+2}-1) \gam_1 +\frac{1}{2} f_{2k} (\be_3 +\gam_3), \\
\gam_3 \be_3&=  f_0 \al_0 + f_1 \al_2 + \cdots f_{2k+1}\al_{4k+2}  \\
& \hspace{2cm} +\frac{1}{2}( f_{2k+2}-1) \be_1 +\frac{1}{2}( f_{2k+2}+1) \gam_1 +\frac{1}{2} f_{2k} (\be_3 +\gam_3), \\
\be_3^2 &= \gam_3^2 = g_0 \al_0 + g_1 \al_2 + \cdots + g_{2k+1} \al_{4k+2}  \\
& \hspace{2cm} + \frac{1}{2} g_{2k+2} (\be_1 +\gam_1) + \frac{1}{2} g_{2k} (\be_3+ \gam_3)  \\
 \end{align*}
All the coefficients of irreducible elements are non-negative integers. 
\end{theorem}
{\bf Proof}. \\
It remains to show that  $f_{2k+2}$ is odd, $ f_{2k}$ is even. In the proof of Theorem \ref{b3g3fusionEVEN}, it has been already proved that $g_{2j}$ is even for any $j$. 

Since $k$ is odd, 
$c_{2k+2} =1$ (mod $4$). Thus
$f_{2k+2}-1= \frac{1}{2} (c_{2k+2}-1)$ is even, i.e. $f_{2k+2}$ is odd.  Likewise $c_{2k}=0$ (mod $4$), thus $f_{2k} =\frac{1}{2} c_{2k}$ is even.  \qed

Thus far we determined that $\be_1$ and $\gam_1$ are self-conjugate, and  computed full irreducible decomposition of $\be_3$, $\gam_3$, in particular ${{\ovl{\be_3}}}= \gam_3$. This determines the rest of the fusion rule. Note that the conjugate map $\pi$ on $\Z V_{11}$ is now determined. 

First, for $\al_{2j}$, $j=0, 1,  \ldots 2k+1$, the right and left multiplication of $\ \al_{2j}$ on any other object from $V_{11}$ is represented by the matrices $Q_j (\DD)$ and $Q_j(\pi \DD \pi)$ respectively. 
\begin{claim}
\dlabel{als}
The entries of the matrices $R_i(\D)$ for $i=0, 1, \ldots 4k+3$ are non-negative integers. In particular, 
the entries of the matrices $Q_j (\DD)$ for $j=0, 1, \ldots 2k+1$ are non-negative integers. 
\end{claim}
\noindent
{\bf Proof.} \\
Immediate from the result in \cite{HW}, which states that, when $\D$ is an adjacency matrix of a graph with norm greater than $2$, then $R_i(\D)$ has non-negative integer entries for any $i$. 
\qed

It remains to determine the decomposition of tensor product of $\be_1$, $\gam_1$ with themselves and $\be_3, \gam_3. $

Since by the graph $\be_1=\be_3 \al_2$, $\gam_1= \gam_3 \al_2$, the fusion among $\be_3$ and  $\gam_3$ together with the fusion of $\al_2$ with all the objects determine $\be_3 \be_1$, $\gam_3 \gam_1$, $\be_3 \gam_1$, $\gam_3 \be_1$ by imposing associativity.  Taking the conjugate, we obtain 
$\be_1 \be_3$, $\gam_1 \gam_3$, $\be_1 \gam_3$, $\gam_1 \be_3$ as well. $\be_1^2 =\be_1 \gam_3 \al_2$, $\gam_1^2 =\gam_1 \gam_3 \al_2$, $\be_1 \gam_1=\be_1 \gam_3 \al_2$, $\gam_1 \be_1= \gam_1 \be_3 \al_2$ are thus all determined. Since there is no division, subtraction of objects are involved in the process of determining each desired fusion rule,  the coefficients are all non-negative integers. 
\subsection{Fusion rules on ${}_N \cX_N \times {}_N \cX_M$}
 \dlabel{fusionNNNM}
We identify $ {}_N \cX_N$ with $V_{11}$,  $ {}_N \cX_M$ with $V_{12}$. From Claim \ref{als}, $\al_i  Y$ for $i$ even and any $Y \in V_{12}$ are determined, so are $X  \al_i$ for $X \in V_{11}$ and $i$ odd. Thus it remains to obtain $\be_i   Y$ and $\gam_i  Y$, where $i=1, 3$, $Y=\be_2$ or $\gam_2$. They are easily determined, since $\be_2 =\be_3   \al_1$, $\gam_2=\gam_3   \al_1$, and the fusion among $\be_i$, $\gam_j$, $i, j=1, 3$ are already determined. Here we imposed associativity again. Since the fusion coefficients among $\be_i$'s and $\gam_j$'s are non-negative integers and product of $  \al_1$ from the right gives fusion with non-negative integers, the fusion coefficients of $\be_i   Y$ and $\gam_i Y$ are non-negative integers as well. 
 %
 %
\subsection{Fusion rules on ${}_N \cX_M \times {}_M\cX_N$}
 \dlabel{NMMN}
 Let $X \in {}_N \cX_M$. Then for $j$ odd, 
 $$X  \bal_j=R_j(\D) X.$$
From Claim \ref{als}, $R_j(\D) X$ is a linear combination of the objects in ${}_N \cX_N$ with non-negative integer coefficients.
 It remains to show that $\be_2   \bbe_2$, $\be_2   \bga_2$, $\gam_2   \bbe_2$ and $\gam_2   \bga_2$ also have this property. It is immediate, since $\bbe_2 = \bal_1 \bbe_3$, $ \bga_2=\bal_1   \bga_3$, $\be_2  \bal=\be_1 + \be_3$, $\gam_2   \bal =\gam_1 + \gam_3$, and all the  fusion rules involved have decompositions into simple objects with $\Z_{\geq 0}$-coefficients. 
\subsection{Fusion rules on ${}_M \cX_M \times {}_M \cX_M $ and ${}_M \cX_M \times {}_M \cX_N$}
 Recall that we have identification  ${}_M \cX_M =V_{22}$ and $ {}_M \cX_N=V_{21}$. 
 Let $\D'$ be the adjacency matrix for $\Gamma'$. Then the fusion rules of the tensor products of $\al'_j$'s for $j=0, 2, \ldots, n-1$, as well as $\oal_k$'s for $k=1, 3, \ldots, n-1$ with any objects in $V_{21} \sqcup V_{22}$ are given by the matrices $R_l(\D')$, where $l=0, 1, \ldots, n$. Similarly to Claim \ref{als}, the entries of $R_l(\D')$ are all non-negative integers. Furthermore, using Frobenius reciprocity,  this also takes care of the coefficients of $\al'_j$'s and $\oal_k$'s in the tensor product of two bimodules. 
\subsection{Fusion rules on  ${}_M \cX_M \times {}_M \cX_M $} 
 The remaining issue is to determine the fusion rule among $f$ and $g$. By observing the Perron-Frobenius weights,  ${\ov f} =f$, ${\ov g}=g$. Since for $j$ even, all the $\al'_j$'s are self-conjugate as well, we have $fg=gf$. 
 \begin{theorem}
 \dlabel{fgfusion}
 \begal
 <f^2, f>&= d_{2k-1}, \;
 <f g, f>=d_{2k}, \\
 <fg, g>&=d_{2k+1}, \;
 <g^2, g>= d_{2k+2}, 
 \end{align*}
 where $d_k$'s are as in the proof of Claim \ref{D}, namely defined by $d_j=d_{j-1}+d_{j-2}+d_{j-3}$, $d_{-1}=0$, $d_0=d_1=1$. 
 \end{theorem}
 \begin{lemma}
 \dlabel{differences}
 \begal
 <f^2, f> - <fg, g>&=d_{2k-1}-d_{2k+1}, \\
 <fg, f>-<g^2, g>&=d_{2k}-d_{2k+2}, \\
 <fg, g>-<g^2, g>&=d_{2k+1}-d_{2k+2}. 
 \end{align*}
 \end{lemma}
 \noindent
 Proof of Lemma \ref{differences}
We use the similar strategy as in Claim \ref{D}. Let $G'$ be the adjacency matrix for $(V_{22}, V_{21})$ corresponding to the graph $\Gamma_k'$, and let 
$$\D' :=
\left(\begin{array}{cc}0 & G'   \\ G'^t & 0 \end{array}\right).
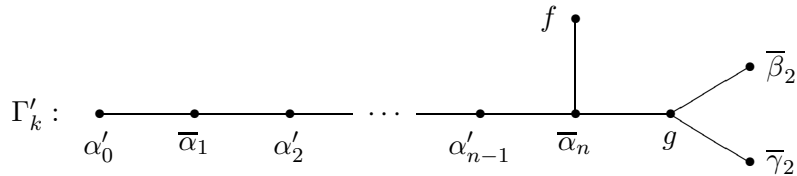
$$ 
 $$\xy 0;/r1pc/:
 (-1, -5)*{\Gamma_k':};
 (1,-5);(9,-5) **\dir{-};
 (11,-5);(19,-5) **\dir{-};
 (1,-5)*\dir{*}, *++!U{\al'_0}; 
 (4,-5)*\dir{*}, *++!U{{\ovl \al_1}}; 
 (7,-5)*\dir{*}, *++!U{\al'_2}; 
 (10,-5)*{\cdots};
 (13,-5)*\dir{*}, *++!U{\al'_{n-1}}; 
 (16,-5)*\dir{*}, *++!U{{\ovl \al_n}};
  (19,-5)*\dir{*}, *++!U{g};
  (16,-5);(16,-2) **\dir{-};
  (16,-2)*\dir{*}, *++!R{f};
   (19,-5);(21.5,-3.5) **\dir{-};
     (19,-5);(21.5,-6.5) **\dir{-};
     (21.5,-3.5)*\dir{*}, *++!L{{ {\ov \beta_2}}};
         (21.5,-6.5)*\dir{*}, *++!L{{ {\ov \gam_2}}};
\endxy,$$
Observe 
\begal
R_0(\D')(g-f)&=(g-f), \\
R_1(\D')(g-f)&={\ov \gamma_2} +{\ov \beta_2}, \\
R_2(\D')(g-f)&= g+f, \\
R_3(\D')(g-f)& = 2{  \al'}_n, \\
R_4(\D') (g-f)& =2 \al'_{n-1} +f+g, \\
\end{align*}
where $\al'_{j}={\ov \al}_j$ for $j$ odd. Then we have
\begin{gather*}
S_j(\D')(g-f)=2\al'_{n-j+3} 
\end{gather*}
 for $j=3, 4, \ldots, n+3$, where the polynomials $S_j$'s are as defined in the proof of Claim \ref{D}. 
On the other hand, 
\begeq
g+f&=&R_{n+1}(\DD')\al'_0 \\
&=&R_{4k+4}(\DD')\al'_0=Q_{2k+2}(\bal_1 \al_1).
\eneq
  Using Lemma \ref{SRrelations}, 
\begal
&(g+f)(g-f) \\
&= (d_0 S_{2(2k+2) }+d_1 S_{2(2k+1)} + \cdots + d_{2k}S_4 +d_{2k+1} R_2 + d_{2k-1} R_0)(\D') (g-f)\\
&= ( \mbox{linear combination of }  \al'_*\mbox{'s} ) + d_{2k+1} (g+f) +d_{2k-1}(g-f) \\
&=( \mbox{linear combination of }  \al'_*\mbox{'s}) + (d_{2k+1}+d_{2k-1})g + (d_{2k+1}-d_{2k-1})f. 
\end{align*}
Therefore we have 
\begin{align*}
<(g-f)(g+f), g>&=<g^2, g> - <f^2, g> \\ & \quad =d_{2k+1}+d_{2k-1}=d_{2k+2}-d_{2k}, \\
<(g-f)(g+f), f>&=<g^2, f> - <f^2, f> =d_{2k+1}-d_{2k-1}. \qquad (\flat 1)
\end{align*}
We obtain further information by investigating $R_2(\D')(g+f)(g-f)$. Note that $R_2(\D')(g+f)=2 \al'_{n-1}+f + 3g.$ Therefore
\begal
&R_2(\D')(g+f)(g-f) \\
&= (2 \al'_{n-1}+f + 3g)(g-f) \\
&= 2\al'_{n-1}(g-f) + 3g^2 -f^2 -2fg \\
&= (  \al'_*\mbox{'s})  + 2(d_{2k}(g+f)+d_{2k-2}(g-f)) + 3g^2 -f^2 -2fg \\
&= ( \al'_*\mbox{'s})  + 2 ( d_{2k}+d_{2k-2}) g +2(d_{2k}- d_{2k-2}) f + 3g^2 -f^2 -2fg \; \; \; (\sharp1) 
\end{align*}
On the other hand, 
\begal
&R_2(\D')(g+f)(g-f) \\
&=R_2(\D') (2 (d_0 \al'_2 +d_1 \al'_4 +  \cdots + d_{2k} \al'_{4k+2}) )+  (d_{2k+1}+d_{2k-1}) R_2(\D') g \\
& \quad + (d_{2k+1}-d_{2k-1}) R_2(\D') f\\
&= ( \al'_*\mbox{'s})  + 2 d_{2k} (f+g) +  (d_{2k+1}+d_{2k-1}) (\al'_{n-1} +f + 2g) \\
& \quad + (d_{2k+1}-d_{2k-1})  (\al'_{n-1} +g )  \\
&=  ( \al'_*\mbox{'s}) + (2 d_{2k} + d_{2k+1} +d_{2k-1})f  + (2 d_{2k} + 3 d_{2k+1} + d_{2k-1}) g.  \; \; \; (\sharp2) 
\end{align*} 
 Comparing $(\sharp1 )$ and $(\sharp 2)$ we obtain
 \begal
 \dlabel{newtwo}
 3 <g^2, g> -<f^2, g> -2 < fg, g> &=3 d_{2k+1} + d_{2k-1}-2d_{2k-2}, \\
 3<g^2, f>-<f^2, f>-2<fg, f> &= d_{2k+1} + d_{2k-1} +2 d_{2k-2} \qquad (\flat2)
 \end{align*}
 Combining the equations $(\flat1)$ and $(\flat2)$ we obtain the statement of the Lemma. Note that we use Frobenius reciprocity  such as $<fg, f>=<f^2, g>$ etc.  \qed
 \begin{lemma}
 $$<g^2, g>=d_{2k+2}, $$
 which implies, together with Lemma \ref{differences}, Theorem \ref{fgfusion}.
 \end{lemma}
 \noindent
 {\bf Proof}
 Since $g=\bbe_2 \al_1 =\bga_2 \al_1$, 
$$
 2g = (\bbe_2+\bga_2) \al_1 ={\ov{(\be_3 + \gam_3) \al_1 }} \al_1 =\bal_1 (\be_3 + \gam_3) \al_1.
 $$
 Also note $\bga_2 = {\ov \gam_3 \al_1}= \bal_1 \be_3$. 
Therefore
\begal
4<g^2, g> & = <\bal_1 (\be_3 + \gam_3) \al_1 \bal_1 (\be_3 + \gam_3) \al_1, \bal_1 \be_3 \al_1> \\
&= <\al_1 \bal_1  (\be_3 + \gam_3) \al_1 \bal_1 (\be_3 + \gam_3) \al_1 \bal_1, \be_3> \\
&=< (\be_3 + \gam_3)^2 (\al_1 \bal_1 )^3, \be_3>\\
&=< (\be_3 + \gam_3)^2, \be_3  (\al_1 \bal_1 )^3>, 
\end{align*}
where we used $\al_1 \bal_1 (\be_3 +\gam_3)= \be_1+\be_3 + \gam_1 + \gam_3 ={\ov{\be_1+\be_3 + \gam_1 + \gam_3} }=
\ov{ (\be_3 +\gam_3) } \al_1 \bal_1 = (\be_3 +\gam_3) \al_1 \bal_1 $. 
By computation using the graph $\Gamma_k$, one obtains
 $$\be_3  (\al_1 \bal_1 )^3 =5 \be_3 +10 \be_1 + 6 \al_{n-1} + 6 \gam_1 + \al_{n-3} + \gam_3.$$
The formula for $ (\be_3 + \gam_3)^2$ is given in Claim \ref{D}. Using it we obtain 
 \begal
& < (\be_3 + \gam_3)^2, \be_3  (\al_1 \bal_1 )^3> \\
 &=8 c_{2k} + 12 c_{2k+1} + 16 c_{2k+2} \\
 &= 4c_{2k+1} + 8 c_{2k+2}  + 8 c_{2k+3} \\
 &= 4 c_{2k+2}  + 4 c_{2k+3} +4 c_{2k+4}= 4c_{2k+5}=4 d_{2k+2}. 
 \end{align*}
 Therefore $<g^2, g>=d_{2k+2}$. 
\subsection{Fusion rules on  ${}_M \cX_M \times {}_M \cX_N$}
\dlabel{fusionMMMN}
 The remaining problem is to determine the fusion rule on $\{f, g\} \times \{ \bbe_2, \bga_2\}$. 
 \begal
 <f \bbe_2, \bbe_2> & =<f, \bbe_2 \be_2> =< f, \oal_1 \be_3^2 \al_1> \\
 &=<\al_1 f \oal_1, \be_3^2 > = <\al_n \oal_1,\be_3^2 >\\
  &=<\be_3^2, \be_1>+<\be_3^2, \gam_1>+ <\be_3^2, \al_{n-1}>.
 \end{align*}
 Using Theorems \ref{b3g3fusionEVEN} and  \ref{b3g3fusionODD} 
$$ <f \bbe_2, \bbe_2>   =  g_{2k+2}+g_{2k+1}. $$
  Both values are non-negative integers. Similarly we obtain
  \begin{align*}
   <f \bbe_2, \bga_2>&=<f \bga_2, \bbe_2> =  f_{2k+2}+f_{2k+1}, \\
   <f \bga_2, \bga_2> &=     g_{2k+2} + g_{2k+1}.
   \end{align*}
   \begal
   <g \bbe_2, \bbe_2>&=<\bbe_2 \al_1 \bbe_2, \bbe_2> =<\oal_1 \bbe_3 \al_1 \oal_1 \bbe_3, \oal_1 \bbe_3> \\
   &=<\al_1 \oal_1 \gam_3 \al_1 \oal_1, \gam_3 \be_3 >=<\ov{\ov{(\gam_1 +\gam_3)} \al_1 \oal_1}, \gam_3\be_3>. 
   \end{align*}
   \begal
   \ov{(\gam_1 +\gam_3)} \al_1 \oal_1& =(\gam_1+\be_3) \al_1 \oal_1 \\
   &=(\al_{n-1} +\be_1 +2\gam_1 + \gam_3) +\be_1 +\be_3\\
   &=\al_{n-1} +2 (\be_1 +\gam_1) + \gam_3 +\be_3 \\
   &=\ov{\al_{n-1} +2 (\be_1 +\gam_1) + \gam_3 +\be_3}.
   \end{align*}
   Thus, using Theorems \ref{b3g3fusionEVEN} and  \ref{b3g3fusionODD}   we obtain
   \begal
   <g \bbe_2, \bbe_2> =\begin{cases}
   f_{2k+1} +2f_{2k+2} + f_{2k} -1 \; \mbox{if } k \mbox{ even} \\
    f_{2k+1} +2f_{2k+2} + f_{2k}  \; \mbox{if } k \mbox{ odd}
    \end{cases}
    \end{align*}
    Similarly, 
    \begal
   <g \bbe_2, \bga_2>& =<g\bga_2, \bbe_2>  \\
  &  =\begin{cases}
   g_{2k+1} +2g_{2k+2} + g_{2k} +2 \; \mbox{if } k \mbox{ even} \\
    g_{2k+1} +2g_{2k+2} + g_{2k}  \; \mbox{if } k \mbox{ odd}, 
    \end{cases} \\
     <g \bga_2, \bga_2> & = <g \bbe_2, \bbe_2>. 
    \end{align*}
    %

 %
 %
 \section{Existence of the fusion algebra}
 \dlabel{existense}
 Let $k \in \N_0$, and put $n=4k+3$ as before. In this section we will reserve the symbols 
 $$(\al_j)_{0 \leq k \leq n}, \; (\be_j)_{1 \leq j \leq 3}, \; (\gamma_j)_{1 \leq j \leq 3} $$
 for elements in a certain bi-graded $\Z$-algebra $\A$ which we define later. Therefore we relabel the vertices of the graph $\Gamma_k$ in the following way: 
 $$
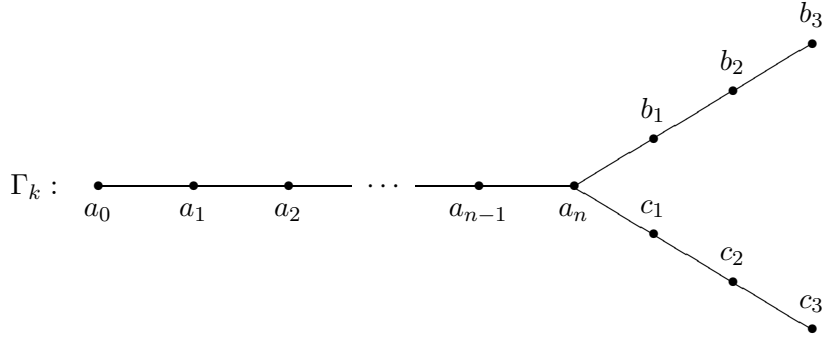
\xy 0;/r1pc/:
  (-1, -5)*{\Gamma_k:};
 (1,-5);(9,-5) **\dir{-};
 (11,-5);(16,-5) **\dir{-};
 (1,-5)*\dir{*}, *++!U{a_0}; 
 (4,-5)*\dir{*}, *++!U{a_1}; 
 (7,-5)*\dir{*}, *++!U{a_2}; 
 (10,-5)*{\cdots};
 (13,-5)*\dir{*}, *++!U{a_{n-1}}; 
 (16,-5)*\dir{*}, *++!U{a_n};
  (16,-5);(23.5,-0.5) **\dir{-};
    (16,-5);(23.5,-9.5) **\dir{-};
     (18.5,-3.5)*\dir{*}, *++!D{b_1}; (21,-2)*\dir{*}, *++!D{b_2}; (23.5,-0.5)*\dir{*}, *++!D{b_3}; 
       (18.5,-6.5)*\dir{*}, *++!D{c_1}; (21,-8)*\dir{*}, *++!D{c_2}; (23.5,-9.5)*\dir{*}, *++!D{c_3}; 
\endxy$$
As in Section \ref{fusioneven}, we let $G$ be the adjacency matrix for $(\Gamma_k^{\rm even}, \Gamma_k^{\rm odd})$, where 
\begeq
\Gamma_k^{\rm even}&=& \{a_0, a_2, \dots, a_{n-1}, b_1, c_1, b_3, c_3\}, \\
\Gamma_k^{\rm odd}&=& \{a_1, a_3, \dots, a_{n}, b_2, c_2\}, \\
\eneq
  we set $\DD=G G^t$, and 
  $${\Delta} :=
\left(\begin{array}{cc}0 & G   \\ G^t & 0 \end{array}\right). 
$$
We set $(q_k)_{k=0}^\infty$ be the sequence of polynomials defined 
 \begeq
q_0(t) &=&t^2-5t+3 \\
q_1(t)&=&   (t-1)(t^3-8t^2+17t-5)  \\
q_k(t)&=&(t^2-4t+2)q_{k-1}(t)-q_{k-2}(t), \; k \geq 2 \\
\eneq
as in Section \ref{fusioneven}. 
Then the characteristic polynomial for $\DD$ is 
$$\chi_k(t)=t^2 (t-2)^2 q_k(t).$$
(cf. Section \ref{fusioneven}). Moreover $q_k(t)$ is a polynomial of degree $2k+2$ with $2k+2$ distinct roots, because by \cite{AY}, either $q_k(t)$ or $q_k(t)/(t-1)$ is an irreducible polynomial. From the recursion formula for the $q_k$-polynomials, one obtains
\begeq
q_k(0)&=&2k+3 \\
q_k(2)&=&(-1)^{k+1} (2k+3), 
\eneq
In particular, $0$ and $2$ are not roots of $q_k$. Let $k \in \N_0$ be now fixed. From the above, we knot that $\chi_k(t)$ has exactly $2k+4$ distinct roots $(t_j)_{k=1}^{2k+4}$, where $t_1 =0$, $t_2=2$, and $t_3, \dots, t_{2k+4}$ are the roots of $q_k(t)$. Since $\DD =G G^t$ is a positive operator, $t_j \geq=0$ for $1 \leq j \leq 2k+4.$ 
\begin{lemma}
\dlabel{uffehand3.1}
Let $E_j$ be the orthogonal projection on the eigenspace of $\DD$ corresponding to the eigenvalue $t_j$ ($1 \leq j \leq 2k+4$) and put 
$$ \mu_j =< E_j a_0, a_0>, $$
where $< \bullet, \bullet>$ is the inner product in $l^2(\Gamma_k^{\rm even})$. 
Then 
\begin{itemize}
\item[(a)] $\Sigma_{j=1}^{2k+4} \mu_j=1$, \\
\item[(b)] $\mu_j >0$ for $1 \leq j \leq 2k+4$, \\
\item[(c)] $\mu_1=\mu_2=\frac{1}{2k+3}$. 
\end{itemize}
\end{lemma}
\noindent
{\bf Proof}. \\
(a): Since $\DD$ is a symmetric matrix, $\Sigma_{j=1}^{2k+4} E_j=I$, which proves (a). \\
(b): From Section \ref{fusioneven}, we have 
\begeq
Q_j(\DD)a_0 &=& R_{2j}(\Delta)a_0=a_{2j} \; (0 \leq j \leq 2k+1), \\
Q_{2k+2}(\DD) a_0&=& R_{4k+4}(\Delta)a_0=b_1 +c_1, \\
Q_{2k+3}(\DD) a_0&=& R_{4k+6}(\Delta) a_0 =b_1 +c_1 +b_3 +c_3. 
\eneq
Since $ \{a_0, a_2, \dots, a_{4k+2}, b_1+c_1, b_1 +c_1 +b_3 +c_3\}$ is a set of $2k+4$ linear independent vectors in $l^2(\Ge)$, and since $(Q_j)_{0 \leq j \leq 2k+3}$ spans the set of polynomials of degree less or equal to $2k+3$, we have
$$ P(\DD)a_0 \not= 0$$
for every non-zero polynomial $P \in \R[x]$ with $\deg(P) \leq 2k+3$. On the other hand, since $\DD$ is diagonalisable with with eigenvalues $(t_j)_{j=1}^{2k+4}$, we have 
$$E_j =P_j(\DD),$$
 where 
 $$P_j(t)= \Pi_{i \neq j} \frac{t-t_i}{t_j-t_i}, \; (t \in \R)$$
 is a polynomial of degree $2k+3$. Hence
 $$\mu_j=<E_k a_0, a_0>= ||E_j a_0||^2 >0 \; (1 \leq j \leq 2k+4). $$
 \noindent
 (c): From Section \ref{fusioneven}, we have 
 \begeq
 \rg(E_1)&=&E(\DD, 0)=\spa\{y_1, y_2\}, \\
 \rg(E_2)&=&E(\DD, 2)=\spa\{x_1, x_2\}, 
 \eneq
where
\begin{align*}
x_1& :=2(a_0 +a_2)-2 (a_4 +a_6)+ \cdots + (-1)^k 2(a_{4k} +a_{4k+2})  \\
& + (-1)^{k+1} (b_1 +c_1 +b_3 +c_3) \\
x_2& :=  (b_1-c_1) +(b_3-c_3) \\
y_1&:= 2a_0 -2 a_2 + \cdots + 2 a_{4k}-2a_{4k+2} + (b_1+c_1)  - (b_3 +c_3)\\
y_2 &:=  (b_1 -c_1) -(b_3 -c_3)
\end{align*}
Since $y_1 \perp y_2$ and $y_2 \perp a_0$, we get 
$$\mu_1=< E_1 a_0, a_0>=\frac{|<y_1, a_0>|^2}{||y_1||^2}=\frac{1}{2k+3}$$
and similarly
$$\mu_2=<E_2 a_0, a_0>=\frac{|<x_1, a_0>|^2}{||x_1||^2}=\frac{1}{2k+3.}$$
\qed
\begin{cor}
\dlabel{uffehand3.2}
Let $(e_{ij})_{i,j=1}^{2k+4}$ be the matrix units of $M_{2k+4}(\R)$. Put 
\begeq
\B&=&\spa_\R\{ e_{11}, e_{12}, e_{21}, e_{22}, e_{33}, e_{44}, \dots, e_{2k+4, 2k+4}\} \\
&\cong& M_2(\R) \oplus l^\infty (\{3, 4, \dots, 2k+4\}, \R)
\eneq
Then $\B$ is a finite dimensional real $C^*$-algebra and $\mu: \B \to \R$ given by 
$$\mu(b):= \sum_{j=1}^{2k+4} \mu_j b_{jj}, \; b=(b_{ij})_{i,j=1}^{2k+4} \in \B$$
is a faithful trace state on $\B$. 
\end{cor}
{\bf Proof}
it is clear from (a) and (b) in Lemma \ref{uffehand3.1} that $\mu$ is a faithful state on $\B$ and the trace property
$$\mu(bc)=\mu(cb), \; b, c \in \B$$
follows from (c) in Lemma \ref{uffehand3.1}.
\qed
\begin{lemma}
\dlabel{uffehand3.3}
Let $k \in \N_0$ be fixed and let $\mu:\B \to \R$ be the trace defined above, and put 
$$A:=\diag (0, \sqrt{2}, \sqrt{t_3}, \dots \sqrt{t_{2k+4}})),$$
 where $t_3, \dots, t_{2k+4}$ are the roots of $q_k$. Then 
 \begin{itemize}
 \item[(a)] For every even polynomial $P \in \R[x]$ 
 $$\mu(P(A))=<P(\Delta)a_0, a_0>.$$
 \item[(b)]
 Let $P, Q \in \R[x]$ be two polynomials, which are either both even or both odd. Then 
 $$\mu(P(A) Q(A))=<P(\Delta) a_0, Q(\Delta)a_0>. $$
 \item[(c)]
 Let $n=4k+3$ (as usual), then 
 $$R_{n+4}(A)-R_{n+2}(A)-R_n(A)-R_{n-2}(A)=0.$$
 \end{itemize}
 \end{lemma}
 \noindent
 {\bf Proof} \\
(a): Let  $Q \in \R[x]$ be so that $P(t)=Q(t^2)$. Then 
 $$<P(\Delta) a_0, a_0>=<Q(\DD)a_0, a_0>.$$
 Let $E_j$ denote the spectral projection of $\DD$ corresponding to the eigenvalue $t_j$ ($1 \leq j \leq 2k+4$) as before, where $t_1=0$ and $t_2=2$. Then 
 $$Q(\DD)=\sum_{j=1}^{2k+4} Q(t_j)E_j.$$
 Hence
 \begeq
 <Q(\DD)a_0, a_0>&=& \sum_{j=1}^{2k+4} Q(t_j) <E_j a_0, a_0>\\
 &=& \sum_{j=1}^{2k+4} \mu_j Q(t_j) \\
 &=&\mu(Q(A^2))=\mu(P(A)).
 \eneq
 \noindent
 (b): Under the assumption on $P$ and $Q$, the product $PQ$ is an even polynomial. Hence by (a) we have
 \begeq
 \mu(P(A)Q(A))&=& < P(\Delta)Q(\Delta)a_0, a_0> \\
 &=& <P(\Delta) a_0, Q(\Delta)a_0>.
 \eneq
 \noindent
 (c): Put $P=Q=R_{n+4} -R_{n+2} -R_n -R_{n-2}$, which is an odd polynomial. By (b), 
 $$\mu(P(A)^2)=||P(\Delta)a_0||_2^2.$$
 From the recursive formula for the polynomials $R_j$ one has
 \begeq
 R_{n-2}(\Delta)a_0&=& a_{n-2}, \\
 R_n(\Delta)a_0&=&a_n,\\
 R_{n+2}(\Delta)a_0&=& a_n + b_2 + c_2, \\
 R_{n+4}(\Delta)a_0 &=&a_{n-2}+2 a_n+b_2 +c_2\\
 &=& (R_{n+2}(A)+R_n(A)+R_{n-2}(A))a_0.
 \eneq
 Hence
 $\mu(P(A)^2)=||P(\Delta)a_0||_2^2 =0$, and since $\mu$ is a faithful trace on $\B$, we have $P(A)=0$. \qed
 \begin{remark}
 Since $P=R_{n+4}-R_{n+2}-R_n-R_{n-2}$ is an odd polynomial and $P(A)=0$, we know that $P(t)$ has at least $n +4 =4k+7$ roots
 $$0, \pm \sqrt{2}, \pm \sqrt{t_3}, \dots, \sqrt{t_{2k+4}},$$
 which are exactly the distinct roots of $t(t^2-2)q_k(t^2)$. Since $P$ and $t(t^2-2) q_k(t^2)$ are both monic polynomial of degree $4k+7$, it follows that 
 $$(R_{n+4}-R_{n+2}-R_n-R_{n-2})(t)=t(t^2-2)q_k(t^2).$$
 It is not hard to prove this identity directly by using the recursion formulas for the polynomials $\{q_k\}$'s and $\{R_j\}$'s. 
 \end{remark}
 \begin{definition}
 \dlabel{uffehand3.5}
 Let $k\in \N_0$, $n=4k+3$, and let $(\B, \mu)$ and $A=\diag(\sqrt{t_1}, \sqrt{t_2}, \dots, \sqrt{t_{2k+4}}) \in \B$ be as before. Let $(f_{ij})_{i,j=1}^2$ be the matrix units in $M_2(\R)$, and put
 $$V:=V_{11} \sqcup V_{12} \sqcup V_{21} \sqcup V_{22}, $$
 where $V_{ij} \subset \B \otimes f_{ij}$ ($i,j =1, 2$) are described as below: \\
 \begin{itemize}
 \item[a)] $V_{11}=\{\al_0, \al_2, \al_4, \dots, \al_{4k+2}, \be_1, \gam_1, \be_3, \gam_3\} $, where
 \begeq
 \al_{2j}&=& R_{2j}(A) \otimes f_{11}, \; 0 \leq j \leq 2k+1, \\
 \be_1&=& \frac{1}{2}(R_{n+1}(A) + \sqrt{2k+3} (e_{12} + e_{21})) \otimes f_{11}, \\
  \gam_1&=& \frac{1}{2}(R_{n+1}(A) - \sqrt{2k+3} (e_{12} + e_{21})) \otimes f_{11}, \\
  \be_3&=& \frac{1}{2}((R_{n+3}-R_{n+1}-R_{n-1})(A) + \sqrt{2k+3} (e_{12} - e_{21})) \otimes f_{11}, \\
 \gam_3&=& \frac{1}{2}((R_{n+3}-R_{n+1}-R_{n-1})(A) - \sqrt{2k+3} (e_{12} - e_{21})) \otimes f_{11} 
 \eneq
 \item[b)]
 $V_{12}=\{\al_1, \al_3, \al_5, \dots, \al_{4k+3}, \be_2, \gam_2\} $ where
 \begeq
 \al_{2j+1} &=& R_{2j+1}(A) \otimes f_{12}, \; 0 \leq j \leq 2k+1, \\
 \be_2 &=& \frac{1}{2}((R_{n+2}-R_n)(A) + \sqrt{2 (2k+3)} e_{12})\otimes f_{12}, \\
 \gam_2 &=& \frac{1}{2}((R_{n+2}-R_n)(A) - \sqrt{2 (2k+3)} e_{12})\otimes f_{12},
 \eneq
 \item[c)]
 $V_{21}=\{\bal_1, \bal_3, \bal_5, \dots, \bal_{4k+3}, \bbe_2, \bga_2\} $
 where
 \begeq
 \bal_{2j+1}&=& R_{2j+1}(A) \otimes f_{21}, \; 0\leq j \leq 2k+1, \\
 \bbe_2&=& \frac{1}{2}((R_{n+2}-R_n)(A) + \sqrt{2 (2k+3)} e_{21})\otimes f_{21}, \\
  \bga_2 &=& \frac{1}{2}((R_{n+2}-R_n)(A) - \sqrt{2 (2k+3)} e_{21})\otimes f_{21}, 
  \eneq
  \item[d)]
  $V_{22}=\{\al_0', \al_2', \dots, \al_{4k+2}', f, g\}$ where
  \begeq
  \al_j'&=&R_{2j}(A) \otimes f_{22}, \; 0 \leq j \leq 2k+1, \\
  f&=&\frac{1}{2}(R_{n-1}+2R_{n+1}-R_{n+3})(A) \otimes f_{22}, \\
  g&=&\frac{1}{2}(R_{n+3}-R_{n-1})(A) \otimes f_{22}. 
  \eneq
  \item[e)]
  The conjugation map $V_{12} \to V_{21}$ and $V_{21} \to V_{12}$ is already defined earlier. For $V_{11}, V_{22}$ all the elements are defined to be self-conjugate except $\be_3$ and $\gam_3$ which are defined to be conjugate of each other. Note that for every $X \in V_{ij}$, the conjugate ${\ov X}$ is equal to $X^*$ (or $X^t$, since all the matrices here are real). 
  \item[f)] We will equip $\R V_{ij} \subset \B \otimes f_{ij}$ with inner products given by 
  $$<b \otimes f_{ij}, c \otimes f_{ij} >_\mu := \mu(c^tb)=\mu(bc^t)$$
  for every $b, c \in \R V_{ij}$ ($i,j =1, 2$). 
  \end{itemize} 
  \end{definition}
 \begin{lemma}
 \dlabel{uffehand3.6}
 Let $i, j \in \{1, 2\}$. For $X, Y \in V_{ij}$, 
 $$
 <X, Y>_\mu = 
 \begin{cases}
 1 \; {\rm if} \; X=Y, \\
 0 \; {\rm if} \; X \neq Y.
 \end{cases}
 $$
 \end{lemma}
 \noindent
 {\bf Proof} 
\\
Let $(b,c)_\mu:=\mu(c^t b)=\mu(bc^t), b, c \in \B$ be the inner product on $\B$ given by $\mu$, and put $||b||_\mu (b, b)_\mu^{1/2}, b \in \B$.  \\
\noindent
a): Case $(i, j)=(1,1)$. It suffices to show that 
$$S_1:=\{R_0(A), R_2(A), \dots, R_{n+1}(A), (R_{n+3}-R_{n+1}-R_{n-1})(A), e_{12}+e_{21}, e_{12}-e_{21}\}$$
is an orthogonal set in $\B$ and that 
\begeq
&&||R_{2j}(A)||_\mu^2=1, \; 0 \leq j \leq \frac{n-1}{2}, \\
&&||R_{n+1}(A)||_\mu^2=2, \\
&&||(R_{n+3}-R_{n+1}-R_{n-1})(A)||_\mu^2=2, \\
&&|| e_{12}+e_{21}||_\mu^2 =|| e_{12}-e_{21}||_\mu^2 =\frac{2}{2k+3}.
\eneq
By the definition of $\mu$ in Corollary \ref{uffehand3.2}, it is clear that $e_{12}+e_{21}$ and $e_{12}-e_{21}$ are $\mu$-orthogonal to the remaining matrices in $S_1$, because $R_j(A)$ is a diagonal matrix for all $j \in \N_0$. 
Moreover, by Lemma \ref{uffehand3.1}, 
\begeq
<e_{12}+e_{21}, e_{12}-e_{21}>_\mu&=&\mu(e_{11}-e_{22})=\mu_1-\mu_2=0, \\
||e_{12}+e_{21}||_\mu^2=||e_{12}-e_{21}||_\mu^2&=& \mu (e_{11} + e_{22})=\mu_1+\mu_2 =\frac{2}{2k+3}. 
\eneq
By Lemma \ref{uffehand3.3} (b), the remaining part of the proof in the $V_{11}$-case reduces to show that 
$$T_1:=\{R_0(\Delta)a_0, R_2(\Delta)a_0, \dots, R_{n+1}(\Delta)a_0, (R_{n+3}(\Delta)-R_{n+1}(\Delta)-R_{n-1}(\Delta))a_0\}$$
is an orthogonal set in $l^2(\Gamma_k)$ with 
\begeq
&&||R_{2j}(\Delta)a_0||^2=1, \; 0 \leq j \leq n-1, \\
&&||R_{n+1}(\Delta)a_0||^2=2, \\
&&||(R_{n+3}-R_{n+1}-R_{n-1})(\Delta)a_0||^2=2. 
\eneq
This follows from the fact that 
$$T_1=\{a_0, a_2, \dots, a_{n-1}, b_1+c_1, b_3+c_3\}. $$
\noindent
b) cases $(i, j)=(1, 2)$ and $(i,j)=(2, 1)$. It suffices to show that 
$$S_2:=\{R_1(A), R_3(A), \dots R_n(A), (R_{n+2}-R_n)(A), e_{12}\}$$
is an orthonormal set in $\B$ and that
\begeq
&& ||R_{2j+1}(A)||_\mu^2=1, \; 0 \leq j \leq \frac{n-1}{2}, \\
&&||(R_{n+2}-R_n)(A)||_\mu^2=2, \\
&&||e_{12}||_\mu^2=\frac{1}{2k+3}. 
\eneq
It is easy to check that $e_{12}$ is orthogonal to the remaining elements of $S_2$ and that 
$||e_{12}||_\mu^2 =(2k+3)^{-1}$ by Lemma \ref{uffehand3.3} (b). The remaining statement about the set $S_2$ follow from the fact that 
\begeq
T_2&=&\{R_1(\Delta)a_0, R_3(\Delta)a_0, \dots, R_n(\Delta)a_0, (R_{n+2}-R_n)(\Delta)a_0\} \\
&=&\{a_1, a_3, \dots, a_n, b_2+ c_2\}
\eneq
is an orthonormal set in $l^2(\Gamma_k)$, and that 
\begeq
&&||a_{2j+1}||^2=1, \; 0 \leq j \leq \frac{n-1}{2}\\
&& ||b_2+ c_2||^2=2. 
\eneq
\noindent
c) Case $(i,j)=(2,2)$. The statement follows in this case if we can show that 
$$S_3:=\{R_0(A), R_2(A), \dots, R_{n-1}(A), \frac{1}{2}(R_{n-1}+2R_{n+1}-R_{n+3})(A), 
\frac{1}{2}(R_{n+3}-R_{n-1})(A)\}
$$
is a $\mu$-orthogonal set in $\B$. By Lemma \ref{uffehand3.3} (b)  this reduces to showing that 
$$T_3:=\{a_0, a_2, \dots, a_{n-1}, \frac{1}{2}(b_1+c_1+b_3+c_3), \frac{1}{2}(b_1+c_1-b_3-c_3)\}$$
is an orthogonal set in $l^2(\Gamma_k)$, which is obvious. \qed
 \begin{theorem}
 \dlabel{uffehand3.7}
 Let $V=V_{11} \sqcup V_{12} \sqcup V_{21} \sqcup V_{22}$ as in Definition \ref{uffehand3.5}. Then $\Z V \subset M_2(\B)$ form a fusion ring, with coeficients given by 
 $$N_{X, Y}^Z= < XY, Z>_\mu,$$ 
 where $ \; X \in V_{ij}, Y \in V_{jk}, Z \in V_{ik}, (i,j,k) \in \{1, 2\}^3, $ and with units $\al_0 \in V_{11}$ and $\al_0' \in V_{22}$. Moreover the graph with vertices $V_{11} \sqcup V_{12}$ obtained by right multiplication by $\al=\al_1$ is $\Gamma_k$ and the graph with vertices $V_{21} \sqcup V_{22}$ obtained by right multiplication $\bal$ is $\Gamma'_k$. 
 \end{theorem}
 \noindent
 {\bf Proof}.\\
 Note that by Lemma \ref{uffehand3.6}, $V_{ij}$ is a linear independent set in $\B \otimes f_{ij}$ for all $i, j \in \{1, 2\}.$ Hence 
 $$\dim(\R V_{11})=|V_{11}|=2k+6$$
 and 
 $$\dim(\R V_{12})=\dim(\R V_{21})=\dim(\R V_{22})=2k+4. $$
 This implies that 
 \begeq
 \R V_{11} &=& \B \otimes f_{11}, \\
 \R V_{12}&=& \spa \{e_{12}, e_{22}, e_{33}, \dots, e_{2k+4, 2k+4}\} \otimes f_{12}, \\
  \R V_{21}&=& \spa \{e_{21}, e_{22}, e_{33}, \dots, e_{2k+4, 2k+4}\} \otimes f_{21}, \\
  \R V_{22}&=& \spa \{e_{11}, e_{22}, e_{33}, \dots, e_{2k+4, 2k+4}\} \otimes f_{22}, \\
 \eneq
 because the four inclusions $\subset$ are obvious, and the right hand sides have dimensions $2k+6$ (resp. $2k+4$, $2k+4$, $2k+4$). Therefore
 $$\R V= \R V_{11} \oplus \R V_{12} \oplus \R V_{21} \oplus \R V_{22}$$
 form a bi-graded $\R$-algebra, and the conjugation $X \to {\ov X} $ extends by linearity to all of $\R V$ and it is given by transposition of matrices. Moreover, for $X \in V_{ij}$, $Y \in V_{jk}, (i,j,k \in \{1, 2\})$, we have a unique decomposition 
 $$XY =\sum_{Z \in V_{ik}} N_{X,Y}^Z Z,$$
 where by Lemma \ref{uffehand3.6}
 $$N_{X, Y}^Z= < XY, Z>_\mu \in \R.$$
 The identities
 $$
 N_{X, Y}^Z=N_{Z, {\ov Y}}^X=N_{{\ov X}, Z}^Y=N_{{\ov Z}, X}^{\ov Y}=N_{Y, {\ov Z}}^{\ov X}$$
 is now a simple consequence of the fact that $\mu$ is a trace state on the real $C^*$-algebra $\B$, so in particular 
 \begeq 
 \mu(b)&=&\mu(b^t), \; b \in \B, \\
 \mu(bc)&=& \mu(cb), \; b, c \in \B
 \eneq
 It remains to be proved that $N_{X, Y}^Z \in \N_0$ and  that multiplication from the right by $\al=\al_1$ (resp $\bal$) on $V_{11}$ (resp $V_{22}$) generates the graph $\Gamma_k$ (resp. $\Gamma'_k$).  
 \begin{lemma}
 \dlabel{uffehand3.8}
 Let $\al=\al_1$.  \\
 a) For $X \in V_{11}, Y \in V_{12}$, 
 $$<X \al, Y>_\mu=<X, Y \bal>_\mu \in \N_0, $$
 and $(<X \al, Y>_\mu)_{X \in V_{11}, Y \in V_{12}}$ is the adjacency matrix $G_k$ for $\Gamma_k$. \\
 b) For $X \in V_{22}, Y \in V_{21}$, 
 $$<X \bal, Y>_\mu=<X, Y \al>_\mu \in \N_0,$$
 and $(<X \bal, Y>_\mu)_{X \in V_{22}, Y \in V_{21}}$ is the adjacency matrix $G'_k$ for $\Gamma'_k$. 
 \end{lemma}
 \noindent
 {\bf Proof} \\
 This follows from simple computations using Definition \ref{uffehand3.5}, Lemma \ref{uffehand3.6}, the recursion formula
 $$(\star) \; \; t R_n(t)=R_{n+1}(t)+R_{n-1}(t), n\geq 1$$
 and the identity from Lemma \ref{uffehand3.3}(c)
 $$(\star \star)\;\; R_{n+4}(A)-R_{n+2}(A)-R_n(A)-R_{n-2}(A)=0:$$
 a) It follows immediately from $(\star)$ that for $1 \leq j \leq 2k+1$, 
 $$\al_{2j}\al =\al_{2j+1} +\al_{2j-1}$$
 which shows that $\al_{2j} \in V_{11}$ is connected to $\al_{2j+1}$ and $\al_{2j-1}$ in $V_{12}$ (with simple edges) and not connected to any other $Y \in V_{12}$. To prove that we recover the graph $\Gamma_k$ this way we just have to check that $\al_0 \al =\al_1$, which is obvious, and that $\be_1 \al=\al_n +\be_2$, $\be_3 \al=\be_2$. 
 The last one follows from 
 \begeq
 \be_3\al&=& \frac{1}{2} ( (R_{n+3}-R_{n+1}-R_{n-1})(A)+\sqrt{2k+3}(e_{12}+e_{21}))A)\otimes f_{12} \\
 &=&\frac{1}{2} (R_{n+4}-2R_n -R_{n-2})(A)+\sqrt{2(2k+3)}e_{12})\otimes f_{12}\\
 &=&\frac{1}{2}((R_{n+2}-R_n)(A)+\sqrt{2(2k+3)}e_{12})\otimes f_{12}\\
 &=& \be_2,
\eneq
where we have used $(\star)$ and $(\star \star)$ and the fact that $e_{12}A=\sqrt{2}e_{12}, e_{21} A=0$. The proof of $\be_1 \al=\al_n +\be_2$ is similar. \\
b) To recover the graph $\Gamma_k$ from $V_{22} \sqcup V_{21}$, it suffices to prove that 
\begeq
\al'_0 \bal &=& \bal_1, \\
\al'_{2j}\bal &=& \bal_{2j+1} +\bal_{2j-1} \; (1 \leq j \leq 2k+1)\\
f \bal &=& \bal_n\\
g \bal &=& \bal_n +\bbe_2 + \bga_2
\eneq
The first two are obvious. Let us prove $f\bal =\bal_n$. The formula for $g \bal$ is obtained in the same way
\begeq 
f\bal &=& \frac{1}{2}((R_{n-1}(A)+2R_{n+1}(A)-R_{n+3}(A))A \otimes f_{21}\\
&=& \frac{1}{2} (R_{n-2}+3R_{n} +R_{n+2} -R_{n+4})(A) \otimes f_{21}\\
&=& \frac{1}{2} \cdot 2 R_n(A) \otimes f_{21}\\
&=& \bal_n
\eneq
where we again have used $(\star)$ and $(\star \star)$. 
\begin{lemma}
\dlabel{uffehand3.9}
Put $$\xi:=(\be_1-\gam_1)+(\be_3 -\gam_3).$$ Then 
$$\bxi:=(\be_1-\gam_1)-(\be_3 -\gam_3), $$ 
and 
\begeq
\frac{1}{2} \xi \bxi &=& 2\al_0 -2 \al_2 + \cdots + 2\al_{4k}-2\al_{4k+2} + (\be_1+\gam_1)  - (\be_3 +\gam_3) \\
\frac{1}{2} \bxi \xi &=& 2(\al_0 +\al_2)-2 (\al_4 +\al_6)+ \cdots + (-1)^k 2(\al_{4k} +\al_{4k+2})  \\
&& + (-1)^{k+1} (\be_1 +\gam_1 +\be_3 +\gam_3) 
\eneq
\end{lemma}
 {\bf Proof}. \\
 Clearly $\bxi=(\be_1-\gam_1)-(\be_3 -\gam_3)$. By Lemma \ref{uffehand3.8}, we know that the linear maps
 \begeq
 &&R_\al: \R V_{11} \to \R V_{12}\\
 &&R_\bal: \R V_{12} \to \R V_{11}
 \eneq
 obtained by right multiplication by $\al$ (resp. $\bal$) have the matrices $G^t$ (resp. $G$) expressed with respect to bases $V_{11}$ for $\R V_{11}$ and $V_{11}$ for $\R V_{12}$. Hence 
 $$R_{\al \bal}:=R_\bal R_\al : \R V_{11} \to \R V_{12}$$ 
 has the matrix $\DD=G G^t$ with respect to the basis $V_{11}$ for $\R V_{11}$. We can now argue exactly as in Case 1 of Section \ref{fusioneven} to get 
 \begeq
 \xi \bxi &\in& E(\DD, 0)_{sc}=\R y_1,\\
  \bxi \xi &\in& E(\DD,2)_{sc}=\R x_1, 
  \eneq
  where
  \begeq
  y_1&=& 2\al_0 -2 \al_2 + \cdots + 2\al_{4k}-2\al_{4k+2} + (\be_1+\gam_1)  - (\be_3 +\gam_3) \\
  x_1& =& 2(\al_0 +\al_2)-2 (\al_4 +\al_6)+ \cdots + (-1)^k 2(\al_{4k} +\al_{4k+2})  \\
&& + (-1)^{k+1} (\be_1 +\gam_1 +\be_3 +\gam_3). \\
\eneq
Since $<\xi \bxi, \al_0>_\mu =<\bxi \xi, \al_0>_\mu=<\xi, \xi>_\mu=4$ and $<y_1, \al_0>_\mu=
<x_1, \al_0>_\mu =2$, it follows that $\xi \bxi=2y_1$ and $\bxi \xi=2 x_1$. \qed
\\ \ \\
{\bf End of proof of Theorem \ref{uffehand3.7}}. 
\\
It remains to be proved that $N_{X,Y}^Z \in \N_0$ for all $X \in V_{ij}, Y \in V_{jk}$ and $Z \in V_{ik}, (i,j \in \{1, 2,3\}). $Having established the formulas for $\xi \bxi$ and $\bxi \xi$ in Lemma \ref{uffehand3.8}, the proof of $N_{X,Y}^Z \in \N_0$ can be obtained from Section \ref{unique}. Using that 
$$N_{X,Y}^Z=N_{Z, {\ov Y}}^X=N_{{\ov X}, Z}^Y,$$
one gets that if $X, Y$ or $Z$ is one of the elements $(\al_j)_{0\leq j \leq n}, (\al'_j)_{0\leq j \leq n}$ (where $\al'_{2k+1}=\bal_{2k+1}$), then $N_{X,Y}^Z$ is an entry of the matrix $R_j(\Delta)$ or $R_j(\Delta')$, which by \cite{HW} is a non-negative integer. In the remaining cases, $X, Y, Z$ are compatible and comes from the list
$$\be_1, \gam_1, \be_3, \gam_3, \be_2, \gam_2, \bbe_2, \bga_2, f, g.$$
For $X, Y, Z \in \{\be_1, \gam_1, \be_3, \gam_3\}, $ we have $N_{X,Y}^Z \in \N_0$ by Theorem \ref{b3g3fusionEVEN}, \ref{b3g3fusionODD}, and the remark at the end of Section \ref{fusioneven}. The case $X, Y, Z \in \{f, g\}$ is treated in Theorem \ref{fgfusion} and the remaining cases can easily be reduced to these two cases by using $\be_2=\be_3 \al$ and $\gam_2=\gam_3 \al$ (c.f. Sections \ref{fusionNNNM} and \ref{fusionMMMN}). \qed
\theoremstyle{definition}
\begin{remark}
From Definition \ref{uffehand3.5}, we have 
\begeq
\xi&=& (\be_1 -\gam_1)+(\be_3 -\gam_3)= 2\sqrt{2k+3} e_{12} \otimes f_{11}, \\
\bxi&=& (\be_1 -\gam_1)-(\be_3 -\gam_3)= 2\sqrt{2k+3} e_{21} \otimes f_{11}.
\eneq
Thus 
\begeq
\xi \bxi &=& 4(2k+3) e_{11} \otimes f_{11} \\
\bxi \xi &=& 4(2k+3) e_{22} \otimes f_{11}. 
\eneq
Since 
$A=\diag (0  ,\sqrt{2}, \sqrt{t_3}, \dots, \sqrt{t_{2k+4}})$, where $t_3, \dots, t_{2k+4}$ are the distinct roots of $q_k(t)$, and since $0, 2 \notin \{t_3, \dots, t_{2k+4}\}$, $e_{11}$ and $e_{22}$  are the projections on the eigenspaces for $A$ with eigenvalues $0$ and $2$ respectively. Using $q_k(0)=2k+3$ and $q_k(2)=(-1)^{k+1} (2k+3)$, one gets
\begeq
(2-A^2) q_k(A^2)&=& 2(2k+3) e_{11}\\
A^2 q_k(A^2)&=&(-1)^{k+1} (2k+3) e_{22},
\eneq
because the polynomial $(2-t)q_k(t)$ vanishes at $t=2$ and $t=t_j, 3 \leq j \leq 2k+4$ and has the value $2(2k+3)$ at $t=0$. Similarly $t q_k(t)$ vanishes at $t=0$ and $t=t_j, 3 \leq j \leq 2k+4$ and has the value $(-1)^{k+1} 2 (2k+3)$ at $t=2$. 
Hence the following two identities holds:
\begeq
\xi \bxi&=& 2(2-A^2) q_k(A^2) \otimes f_{11} = 2(1_N -\al \bal) q_k (\al \bal) \\
\bxi \xi &=& (-1)^{k+2} 2 A^2 q_k (A^2) \otimes f_{11} = (-1)^{k+2} 2 \al \bal q_k(\al \bal),\\
\eneq
where $1_N =\al_0$ and $\al=\al_1$. Let $Q_j$ denote as usual the polynomial for which $R_{2j}(t)=Q_j(t^2), t \in \R$. Then by Definition \ref{uffehand3.5}, 
\begeq
\al_{2j} &=& Q_j (\al \bal) \\
\be_1 +\gam_1 &=& Q_{2k+2} (\al \bal)\\
\be_3 +\gam_3 &=& (Q_{2k+3} - Q_{2k+2} -Q_{2k+1})(\al \bal).\\
\eneq
Hence a more direct proof of Lemma \ref{uffehand3.8} can be obtained if the two polynomial identifies (i) and (ii) below holds: Put 
$$r_k(t)=(2-t) q_k (t), \; s_k(t) =(-1)^{k+1} t q_k(t). $$
Then 
\begeq
(i) \;\; r_k &=& (2 Q_0 -2 Q_1 + \cdots + 2Q_{2k}-2Q_{2k+1}) \\
 && +(Q_{2k+1} +2 Q_{2k+2} -Q_{2k+3})\\
(ii) \;\; s_k &=&2(Q_0 +Q_2) - 2 (Q_2 +Q_4) + \cdots + (-1)^k 2(Q_{2k} +Q_{2k+1}) \\
 &&+ (-1)^{k+1} (Q_{2k+3} -Q_{2k+1}). 
\eneq
These two polynomials identities are actually true, and they can be proved by using the recursion formulas for $(q_k)_{k=0}^\infty $ and $(R_j)_{j=0}^\infty$. \qed
\end{remark}

 \thebibliography{999} 
 
 \dbibitem{A2} 
Asaeda, M. (2007).
Galois groups and an obstruction to principal graphs of subfactors.
{\em International Journal of Mathematics}, {\bf 18}, 191--202.

 \dbibitem{AH}
Asaeda, M. and Haagerup, U. (1999).
Exotic subfactors of finite depth with Jones indices
${(5+\sqrt{13})}/{2}$ and ${(5+\sqrt{17})}/{2}$.
{\em Communications in Mathematical Physics},
{\bf 202}, 1--63.

\dbibitem{AY}
Asaeda, M. and Yasuda, S. (2009). 
On Haagerup's list of potential principal graphs of subfactors. 
{\em Communications in Mathematical Physics} {\bf  286}, 1141--1157.

\dbibitem{BMPS}
Bigelow, S., Morrison, S., Peters, E., Snyder, N., (2009). Constructing the extended Haagerup planar algebra. arXiv: 0909.4099 [math.OA].

\dbibitem{ch}
 ErdŽlyi, E. (1953).  Higher Trancsendental Functions, Vol II, McGraw-Hill, section 10.11

\dbibitem{H5}
Haagerup, U. (1994).
Principal graphs of subfactors in the index range 
$4< 3+\sqrt2$. in {\em Subfactors ---
Proceedings of the Taniguchi Symposium, Katata ---},
(ed. H. Araki, et al.),
World Scientific, 1--38.

\dbibitem{HW}
de la Harpe, P. and Wenzl, H. (1987).
Operations sur les rayons spectraux de matrices
symetriques entieres  positives.
{\em Comptes Rendus de l'Academie des Sciences, S\'erie I,
Math\'ematiques},
{\bf 305}, 733--736.
 
 \dbibitem{I1}
Izumi, M. (1991).
Application of fusion rules to
classification of subfactors.
{\em Publications of the RIMS, Kyoto University},
{\bf 27}, 953--994.

\dbibitem{SV}
Sunder, V. S. and Vijayarajan, A. K. (1993).
On the non-occurrence of the Coxeter graphs 
$\beta_{2n+1}$, $E_7$, $D_{2n+1}$ as
principal graphs of an inclusion of II$_1$ factors.
{\em Pacific Journal of Mathematics} {\bf 161}, 185--200.
 
\endthebibliography
 
 \end{document}